\documentclass[review,11pt]{elsarticle}

\makeatletter
\def\ps@pprintTitle{%
 \let\@oddhead\@empty
 \let\@evenhead\@empty
 \def\@oddfoot{}%
 \let\@evenfoot\@oddfoot}
\makeatother

\usepackage{amssymb}

\usepackage{mathtools}
\usepackage{amsmath} %
\usepackage{amsfonts}
\usepackage{latexsym}
\usepackage{siunitx}
\usepackage[outdir=./]{epstopdf}
\usepackage{epsfig}
\usepackage{isomath}
\usepackage{bm, bbm}
\usepackage{xcolor, colortbl}
\usepackage{booktabs}
\usepackage[T1]{fontenc}
\usepackage[utf8]{inputenc}
\usepackage{algorithm, algpseudocode}
\usepackage{setspace}
\usepackage{geometry}
\usepackage{tabularx}
\usepackage{comment}
\definecolor{Gray}{gray}{0.90}
\newcolumntype{a}{>{\columncolor{Gray}}r}
\newcolumntype{L}{>{\raggedright\let\newline\\\arraybackslash\hspace{0pt}}X}
\newcolumntype{R}{>{\raggedleft\let\newline\\\arraybackslash\hspace{0pt}}X}
\usepackage{hyperref}
\usepackage[nameinlink]{cleveref}

\providecommand{\doi}[1]{%
	\begingroup
	\let\bibinfo\@secondoftwo
	\urlstyle{rm}%
	\href{http://dx.doi.org/#1}{%
		doi:\discretionary{}{}{}%
		\nolinkurl{#1}%
	}%
	\endgroup
}

\newcommand{\tensor}[1]{\boldsymbol{\mathbf{#1}}}             %
\renewcommand{\vec}[1]{{\mathbf{#1}}}  %
\newcommand{\Rvec}[1]{\underline{#1}}                %
\newcommand{\Grad}{\operatorname{Grad}}      %

\newcommand{\pcfem}[1][]{%
	\ifthenelse{\equal{#1}{}}{\Rvec{p}_{c}}{\Rvec{p}_{c,{#1}}}%
}
\newcommand{\pfem}[1][]{%
	\ifthenelse{\equal{#1}{}}{\Rvec{p}}{\Rvec{p}_{#1}}%
}
\newcommand{\dd}{\,\mathrm{d}}

\newcommand{\dX}{\,\mathrm{d}\vX}

\newcommand{\dsX}{\,\mathrm{d}s_{\vX}}

\newcommand{\eff}{\mathrm{eff}}
\newcommand{\normalout}{{\vec{n}^{\mathrm{out}}_0}}
\newcommand\aug{\fboxsep=-\fboxrule\!\!\!\vrule\!\!\!}

\newcommand*{\der}{\mathrm{d}}

\newcommand*{\mbF}{\vec{\mathcal{F}}}

\newcommand*{\vu}{\vec{u}}
\newcommand*{\vn}{\vec{n}}
\newcommand*{\vb}{\vec{b}}
\newcommand*{\vv}{\vec{v}}
\newcommand*{\vx}{\vec{x}}
\newcommand*{\vX}{\vec{X}}
\newcommand*{\vZero}{\vec{0}}
\newcommand*{\vf}{\vec{f}}
\newcommand*{\vs}{\vec{s}}

\newcommand{\LV}{lv}  %
\newcommand{\RV}{rv}  %
\newcommand{\LA}{la}  %
\newcommand{\RA}{ra}  %

\definecolor{Federica}{RGB}{0., 145, 240}
\definecolor{Gernot}{rgb}{0.65, 0.04, 0.77}
\definecolor{Christoph}{rgb}{0.94, 0.0, 0.25}
\definecolor{Jordi}{RGB}{0,128,0}

\begin{document}

\begin{frontmatter}

\title{A coupling strategy for a 3D-1D model of the cardiovascular system
to study the effects of pulse wave propagation on cardiac function}
\author[1,2]{Federica Caforio\corref{cor1}}
\author[2]{Christoph M. Augustin}
\author[3,4]{Jordi Alastruey}
\author[2]{Matthias A. F. Gsell}
\author[2]{Gernot Plank}
\address[1]{%
  Institute of Mathematics and Scientific Computing, NAWI Graz,
  University of Graz, Graz, Austria}
\address[2]{%
  Gottfried Schatz Research Center: Division of Biophysics,
  Medical University of Graz, Graz, Austria}
\address[3]{%
Department of Biomedical Engineering, Division of Imaging Sciences and Biomedical Engineering,
King's~College~London, King's~Health~Partners, St.~Thomas'~Hospital, London SE1 7EH, UK}
\address[4]{%
World-Class Research Center ``Digital Biodesign and Personalized Healthcare'', Sechenov University, Moscow, Russia}
\vspace{-3mm} 
\cortext[cor1]{                                                                
  Corresponding author. Email: federica.caforio@uni-graz.at.       
}     
\begin{abstract}
The impact of increased stiffness and pulsatile load on the circulation and their influence on heart performance have been documented not only for cardiovascular events but also for ventricular dysfunctions.
For this reason, computer models of cardiac electromechanics (EM) have to integrate
effects of the circulatory system on heart function to be relevant for clinical applications.
Currently it is not feasible to consider three-dimensional (3D) models of the entire circulation.
Instead, simplified representations of the circulation are used,
ensuring a satisfactory trade-off between accuracy and computational cost.
In this work, we propose a novel and stable strategy to couple a 3D EM model of
the heart to a one-dimensional (1D) model of blood flow in the arterial system.
A personalised coupled 3D-1D model of LV and arterial system is built and used in a numerical benchmark to demonstrate robustness and accuracy of our scheme over a range of time steps. 
Validation of the coupled model is performed by investigating the coupled system's physiological response to variations in the arterial system affecting pulse wave propagation, comprising aortic stiffening, aortic stenosis or bifurcations causing wave reflections.
Our results show that the coupled 3D-1D model is robust, stable and correctly replicates known physiology. 
In comparison with standard coupled 3D-0D models, additional computational costs are negligible,
thus facilitating the use of our coupled 3D-1D model as a key methodology in studies where wave propagation effects are under investigation. 
\end{abstract}
\end{frontmatter}
\section{Introduction}%
\label{sec:introduction}
Computational modelling of cardiac EM is increasingly proposed and pursued as an effective tool for investigating cardiac physiology and numerous pathological conditions~\cite{ augustin2016anatomically, niederer2019computational}.
Computer models have been used for quantitative analysis, since they have the potential to assess valuable diagnostic information and represent a strong predictive tool for analysing the patient-specific response to a given treatment~\cite{diaz2013roadmap}. 
However, this is a very open and challenging field of research, since cardiac
function relies on complex biophysical aspects that are strongly
interlinked~\cite{chabiniok2016multiphysics}.
Indeed, the muscle contraction is stimulated by electrical activation and
strongly interacts with intraventricular blood and the circulatory system
to transport nutrients and clear waste products.
As a consequence, efficient computational tools are required that account for the
complex physiology and multiphysics of the heart and can be practically implemented
and integrated into the clinic~\cite{corral2020digital}.
In this context, numerous mathematical models have been reported in the literature that aim to model the effect of the circulatory system on cardiac function.
These range from simple lumped zero-dimensional (0D)
Windkessel-type models~\cite{westerhof2009arterial,segers2003systemic,segers2008three},
to more refined 1D models~\cite{Formaggia_etal_2003,Reymond_etal_2011, Manganotti2021}.
While 0D models are suitable as global models of afterload, providing a simple way to evaluate global cardiovascular dynamics, they do not consider any effects due to pulse wave transmission and thus are unable to predict important markers such as pulse wave velocity.
In contrast, 1D blood flow models enable an efficient representation of the
effects of pulse wave propagation and reflection in the circulation.
Consequently, 1D models may be preferred over 0D models when local vascular changes
or distributed properties, e.g., vessel branching, tapering, stenoses, and their
effect on central pressure waveforms are studied and when pulse wave
transmission effects are under investigation~\cite{van2011pulse}.
To date, most 3D computational models of the heart proposed in the literature
consider lumped 0D models as boundary conditions to model the
circulation~\cite{sainte2006modeling,kerckhoffs2007coupling,nordsletten2011coupling,augustin2016anatomically,hirschvogel2017monolithic,regazzoni2020cardiac}
for ease of implementation and calibration.
To the best of our knowledge, no studies propose an effective and stable method to
couple a 3D cardiac model to a 1D description of the circulation,
which can provide new insights into the physiology and analyse the effects of
disrupted wave propagation on cardiovascular dynamics.
This is clinically relevant, since the impact of increased stiffness
and pulsatile load on the circulation and cardiac performance have been
clinically documented not only for cardiovascular events,
but also for left and right ventricle dysfunctions~\cite{stevens2012rv}.
In order to address this need, in this methodological paper we present a novel
\emph{in silico} model based on a multidimensional approach, namely a 3D EM model
of cardiac function together with a 1D model of the human circulation. 
The building blocks of the coupled model and the main aspects of model parameterisation are thoroughly described in the text.
The coupling strategy is inspired by recent methods proposed for 3D-0D
models~\cite{sainte2006modeling,gurev2015high} and is based on the solution of
a saddle-point problem for the volume and pressure in the cardiac cavity.
A personalised 3D EM model of the LV~\cite{augustin2016anatomically} coupled to a 1D circulatory model~\cite{Alastruey2012a} is built and employed in numerical test cases to demonstrate robustness and stability of the numerical scheme. 
The physiological response of the coupled model is investigated in different conditions affecting pulse wave transmission, such as aortic stiffening, aortic narrowing and complex vessel networks involving numerous bifurcations.
We show the ability of the coupled 3D-1D model to correctly replicate known physiological behaviours related to pulse wave propagation. 
The computational cost of our coupled 3D-1D model is comparable to that of standard 3D-0D models, suggesting the use of our 3D-1D model in clinical applications where wave transmission effects are investigated.    
\section{Methodology}%
\label{sec:methods}
\subsection{3D Electromechanical cardiac model}%
\label{sec:heart}
Cardiac tissue is modelled as a hyperelastic, anisotropic, nearly incompressible material endowed with a nonlinear stress-strain relationship~\cite{guccione1995finite}.
Stresses associated with active contraction are characterised as orthotropic, and it is assumed that full contractile force acts along the myocyte fibre orientation, whereas 40\% contractile force is directed along the sheet orientation~\cite{Genet2014}.
Active stress generation is represented based on a simplified phenomenological contractile model~\cite{niederer11:_length}.
Active stresses in this model are only length-dependent, and dependence on fibre velocity is neglected. %
A recent reaction-eikonal (R-E) model~\cite{neic17:_efficient} is considered for the generation of the electrical activation sequences that trigger active stress generation in the tissue.
The hybrid R-E model is a combination of a standard reaction-diffusion model (based on the monodomain equation) and an eikonal model.
Spatio-temporal discretisation of all PDEs and the solvers for the resulting systems
of equations are based on the Cardiac Arrhythmia Research Package~\href{https://carpentry.medunigraz.at/carputils/index.html }{(CARPentry)}~\cite{augustin2016anatomically, neic17:_efficient,vigmond08:_solvers}, built upon extensions of the freely available
openCARP EP framework (\url{http://www.opencarp.org}).
Further information on the EM model of cardiac function can be retrieved in~\ref{sec:app_heart}.
For numerical details on the finite element discretisation, we refer the reader to~\ref{sec:fe_varf}.
\subsection{Circulatory system}%
\label{sec:blood_flow}
The circulatory system imposes a load on the heart and, therefore, affects its mechanical activity.
However, the interaction between the heart and the circulatory system is bidirectional, i.e., the outflow of blood from a heart cavity and pressure in the cavity depend on the current state of the circulatory system, and pressure and flow in the cardiovascular system are determined by the current state of the cavity itself.
The full physics of this interaction is most accurately posed as a fluid-structure interaction (FSI) problem where pressure, $p(\vX)$, and flow velocity, $\vv(\vX)$, of the fluid are the coupling variables~\cite{karabelas2018towards,fernandez2007projection,gerbeau2012p}.
Any perturbation in blood flow velocity and pressure changes the state of deformation of the heart and attached vessels. As a result, this change in strain implies a change in stress within the myocardial muscle.
Conversely, any change in strain or compliance of the attached vessel or the heart changes the pressure and flow of the fluid.
Even though such a distributed PDE-based approach may accurately describe this interaction, at the level of the whole circulatory system it is hardly feasible, for computational and structural reasons.
In fact, a multi-beat simulation using a 3D FSI model solution would
be inconvenient for clinical application.
Moreover, the use of 3D models of haemodynamics would require the identification of a much larger number of parameters, that is not easily attainable within clinical constraints.
A reduced-order approximation of the circulatory system is used instead in this work, which based on a 1D model for the circulatory system and is particularly suitable for simulating blood flow in the aorta and the major arteries.
The parameters in the circulatory model are identified and constrained by imaging-based measurements of the heart and blood flow and invasive blood pressure measurements.
In this framework the nonlinear PDE EM model of the heart is coupled to the 1D model of the circulatory system using the lumped hydrostatic pressure $p_\mathrm{cav}$ in the cardiac cavity and the flow $q_\mathrm{cav}$ out of the cavity into the circulatory system as coupling variables.
\subsubsection{1D mathematical model of the blood flow circulation in the arterial system}
The 1D equations of blood flow in the human circulation can be derived from the Navier-Stokes equations after assuming axisymmetric flow in a cylindrical tube with thin wall~\cite{Formaggia_etal_2003}. In this work, we only take into account the arterial circulation and we consider a Voigt-type visco-elastic constitutive law for the arterial wall~\cite{Alastruey2012a}.
A three-element Windkessel lumped parameter model is employed as a terminal boundary condition to simulate the effect of peripheral vessels on pulse wave propagation in larger arteries.
The numerical discretisation of the resulting system of equations relies upon a FEM discretisation (discontinuous Galerkin) in space and finite difference (explicit second-order Adams-Bashforth) in time~\cite{Alastruey2012a}.
The numerical solution of the resulting hyperbolic system of equations is performed with the solver Nektar1D (\url{http://haemod.uk/nektar}). %
We refer the reader to~\ref{sec:oneDmodelling} for more information on 1D blood flow modelling.
The arterial network is connected to the left ventricle (LV) cavity through a model of aortic valve (AV) dynamics, based on Bernoulli's equation~\cite{Mynard_etal_2012}. See~\ref{sec:valves} for further details on the valve dynamics model.
\section{3D-1D coupling}
\label{sec:3D1Dcoupling}
Following the approach for 3D-0D coupling proposed in~\citet{augustin2021computationally} and inspired by \citet{gurev2011models, kerckhoffs2007coupling}, the coupling condition between an EM model of the heart cavities and a reduced-order circulatory model
imposes that the volume change in each cavity (left ventricle LV, right ventricle RV, left atrium LA and right atrium RA)
is balanced with the volume change in the attached circulatory system.
For the sake of generality, we cast in what follows the general framework to
couple a four-chamber heart model with a closed-loop circulatory system.
Let us consider that approximately, at a given time-point, we have constant pressures in each cavity $p_c$,
with $c\in\{\mathrm{LV,RV,LA,RA}\}$.
Then, the following equations are the starting point of the method:
\begin{equation} \label{eq:VolumeEqu}
  V_c^\mathrm{heart}(\vu)-V_c^\mathrm{CS}(p_c)=0,
\end{equation}
where $V_c^\mathrm{heart}(\vu)$ is the cavity volume computed
with the deformation using Eq.~\eqref{eq:volumeOmega},
and $V_c^\mathrm{CS}(p_c)$ is the volume
as predicted by the blood flow model for the intra-cavitary pressure $p_c$. 
We write $\pcfem={[p_{c}]}_{c\in\{\mathrm{LV,RV,LA,RA}\}}$
for the vector of up to $1 \le N_\mathrm{cav} \le 4$ pressure unknowns.
Note, however, that in a purely EM simulation framework the fluid domain is
not explicitly modelled.
Therefore, the cavitary volume $V_c^\mathrm{heart}(\vu)$ is not discretised.
Instead, only the surface enclosing the cavitary volume is known.
Nonetheless, if we assume that the entire surface of the cavitary volume is available, including the faces representing the valves, %
then we can compute for each instant $t$ of the cardiac cycle the enclosed volume $V^\mathrm{heart}_c(\vu)$ from the surface
$\Gamma_c$ by means of the divergence theorem
\begin{equation}  \label{eq:volumeOmega}
  V^\mathrm{heart}_c(\vu)=\frac{1}{3}\int_{\Gamma_c} \vu(t)\cdot\vn\dd \Gamma_c.
\end{equation}
Then, the approach used to evolve the full system of equations of nonlinear elasticity Eq.\eqref{eq:VolumeEqu} together with the coupling condition Eq.\eqref{eq:volumeOmega} is based on the resolution of a saddle-point problem in the variables $(\vu,\pfem)$, corresponding to the displacement and the pressure in the cavity:
\begin{equation}
\begin{split}
  \langle \mathcal{A}_0(\vu),\vv\rangle_{\Omega_0}
  - \langle \mathcal{F}_0(\vu,p_c),\vv\rangle_{\Omega_0} &=\vZero\\
  \langle V_c^\mathrm{heart}(\vu),q\rangle_{\Omega_0}
  - \langle V_c^\mathrm{CS}(p_c), q\rangle_{\Omega_0}  &= 0,
\end{split}
\label{eq:saddlepoint_pb}
\end{equation}
which is valid for all vector fields $\vv$ smooth enough and satisfying the given boundary conditions,
test functions $q$ that are $1$ for the cavity $c$ and $0$ otherwise,
the duality pairing $\langle \cdot ,\cdot \rangle_{\Omega_0}$,
and cavities $c\in\{\mathrm{LV,RV,LA,RA}\}$.
 Note that $\langle A_0(\vu), \vv \rangle_{\Omega_0}$ can be physically interpreted as the rate of internal mechanical work, whereas $\langle \mbF_0(\vu, p_c), \vv \rangle_{\Omega_0}$ takes into account the contribution of pressure loads.
Using a Galerkin discretisation and the Newton method, the problem to be solved at each Newton--Raphson step becomes a block system
in order to find $\delta\Rvec{u}\in\mathbbm{R}^{3N}$ and
$\delta \pfem\in\mathbbm{R}^{N_\mathrm{cav}}$ such that:
\begin{equation}\label{eq:blockCVSystem}
  \begin{pmatrix}
    (\tensor{A}'-\tensor{M}')(\Rvec{u}^k, \Rvec{p}_c^k)  & \tensor{B}'_\mathrm{p}(\Rvec{u}^k) \\
        \tensor{B}'_{\vu}(\Rvec{u}^k) & \tensor{C}'(\Rvec{p}_c^k) \\
 \end{pmatrix}
\begin{pmatrix}
  \delta\Rvec{u}\\ \delta \Rvec{p}_c
\end{pmatrix}
=-
\begin{pmatrix}
  \Rvec{A}(\Rvec{u}^k)-\Rvec{B}_\mathrm{p}(\Rvec{u}^k,\Rvec{p}_c^k) \\
  \Rvec{V}_c^\mathrm{heart}(\Rvec{u}^k) -\Rvec{V}_c^\mathrm{CS}(\Rvec{p}_c^k)
\end{pmatrix},\\
\end{equation}
with the updates
\begin{equation}
    \Rvec{u}^{k+1}   = \Rvec{u}^{k} + \delta\Rvec{u},\quad
    \Rvec{p}_c^{k+1} = \Rvec{p}_c^k + \delta\Rvec{p}_c,
\end{equation}
and the solution vectors $\Rvec{u}^k\in\mathbbm{R}^{3N}$ and
$\Rvec{p}_c^k\in\mathbbm{R}^{N_\mathrm{cav}}$ at the $k$-th Newton
step.
In particular, we can retrieve the expression of $\tensor{C}'(\Rvec{p}_c^k)$ and $\Rvec{V}^\mathrm{CS}(\Rvec{p}_c^k)$ from the equations of the circulatory model, whereas the term $\tensor{B}'_\mathrm{u}(\Rvec{u}^k)$ only depends on the contribution of the EM model of the heart.
In more detail, $\tensor{C}'(\Rvec{p}_c^k)$ is approximated by a discrete derivative (finite difference) of the volume $\operatorname{V}^\mathrm{CS}$ with respect to the cavity pressure.
We refer the reader to~\ref{sec:coupling} and~\ref{sec:FE} for further information on the coupling strategy and finite element formulation, respectively.
\section{Simulation setup}
\label{sec:sim_setup}
\subsection{Numerical framework}
The coupling scheme and a C\texttt{++} circulatory system module have been embedded in CARPentry, based on the software \href{http://haemod.uk/nektar}{Nektar1D}, in order to preserve computational efficiency and strong scalability.
A Schur complement approach (see~\ref{sec:SchurComplement}) is considered to cast
the coupling problem in a pure displacement formulation,
in order to use solver methods already implemented~\cite{augustin2016anatomically}.
A generalised minimal residual method (GMRES) with an
relative error reduction of $\epsilon=10^{-8}$ is employed.
The library \emph{PETSc}~\cite{petsc-user-ref}
and the incorporated solver suite \emph{hypre/BoomerAMG}~\cite{henson2002boomeramg} are used for efficient preconditioning.
A generalised-$\alpha$ scheme is considered for time integration, with spectral radius $\rho_\infty=0$ and damping
parameters $\beta_\mathrm{mass}=\num{0.1}$, $\beta_\mathrm{stiff}=\num{0.1}$, see  \cite{augustin2021computationally} for further information. 
\subsection{Model parameterisation}
\begin{figure}[ht]
	\begin{center}
		\includegraphics[width=8cm]{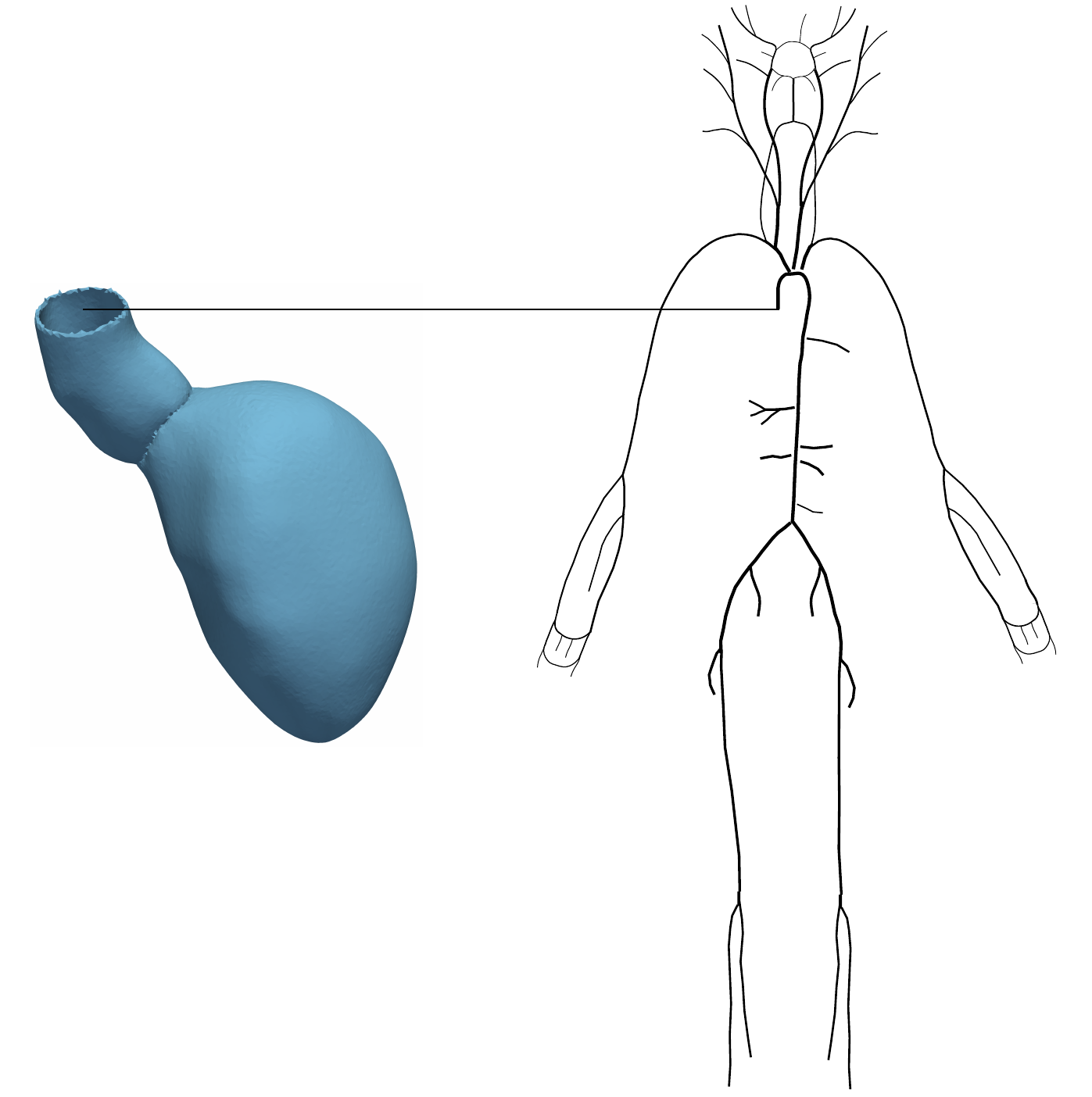}
		\caption{Left: Computational mesh of LV derived from patient-specific image-based clinical data. Right: 1D arterial network.}%
		\label{fig:CoA1}
	\end{center}
\end{figure}

In order to guarantee the stability of the coupled system, the solution of the arterial 1D model is initialised until a periodic solution is reached.
This is achieved by performing 20 cycles of the circulatory model alone prior to coupling with the heart model, with a flow or pressure profile provided as an inlet boundary condition.
For example, a flow or pressure waveform clinically measured at the aortic root can be used as a  boundary condition, if available.
Otherwise, analytical functions can be prescribed~\cite{boileau2015benchmark}.
For the current study, we considered two different choices: an invasively-recorded pressure measurement at the aorta and an analytical function of flow waveform calibrated to match measured peak flow, opening, maximum and closing pressure at the aorta.
The coupled system is solved using a fully converging Newton method with maximal number of steps $k_{\max}=10$ and an absolute $\ell_2$ norm error reduction of the residual of $\epsilon=10^{-6}$.
To illustrate the performance of the coupling method we consider a system consisting of a 3D EM model of the LV and for the 1D arterial model we analyse different geometries, either consisting of a single arterial segment or composed by more complex networks including the major vessels.
To reproduce a realistic case, the LV computational domain is adapted from a patient-specific 3D-whole-heart-MRI scan collected and post-processed in a previous study, see~\Cref{fig:CoA1}. %
The biomechanical parameters of the solid and fluid models under baseline conditions have been calibrated to match a patient-specific set of \emph{in vivo} measurements at specific instants of the cardiac cycle for the same subject, see~\cite{marx2020personalization} for further information. 
The same pressure-volume data is used as for the initialisation of the 3D EM model.
All parameters used for the baseline case are summarised in~\Cref{tab:heart_input_params}.
\section{Results}%
\label{sec:results}
\subsection{Stability assessment of the coupling method}
\label{sec:res_stability}
\begin{figure}[ht]
	\centering
	\begin{center}
		\includegraphics[width=15cm,height=5cm]{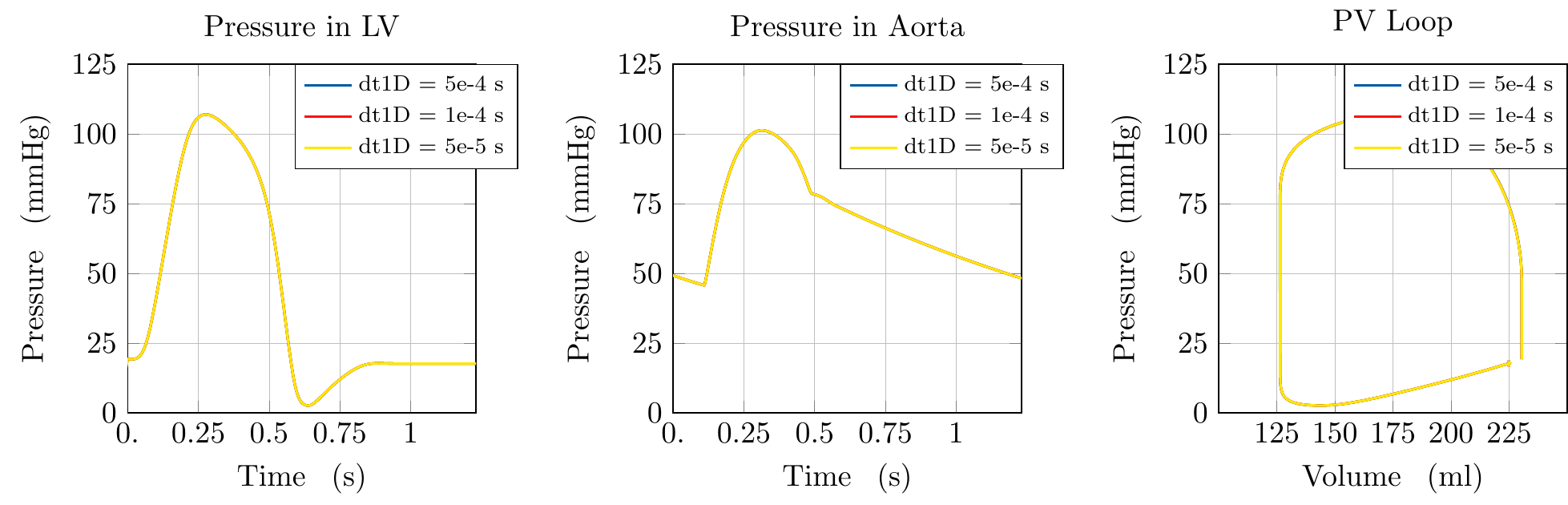}
		\caption{Illustration of model predictions. Comparison of model predictions considering different time resolutions dt1D for the 1D circulation model. Left: Pressure trace in the LV and at the inlet of the aorta. Right: Pressure-volume loop in the LV.}
		\label{fig:res_dt1d}
	\end{center}
\end{figure}

\begin{figure}[ht]
	\centering
	\begin{center}
		\includegraphics[width=15cm,height=5cm]{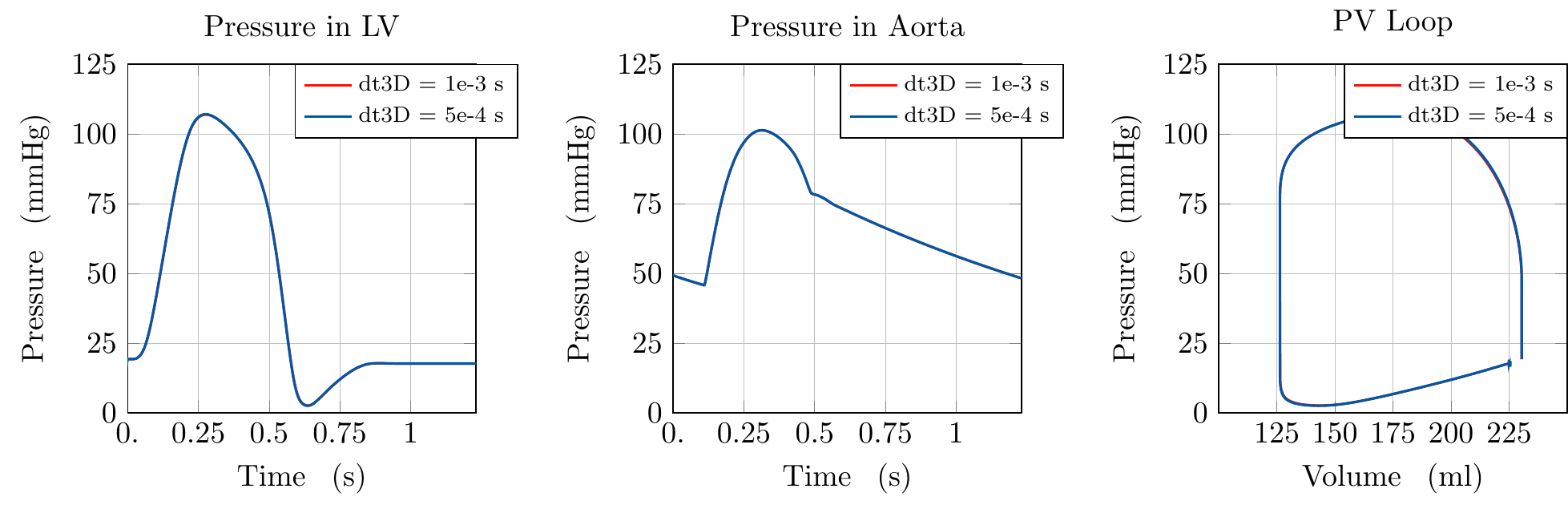}
		\caption{Illustration of model predictions. Comparison of model predictions considering different time resolutions dt3D for the 3D cardiac model. Left: Pressure trace in the LV and at the inlet of the aorta. Right: Pressure-volume loop in the LV.}
		\label{fig:res_dt3d}
	\end{center}
\end{figure}

To analyse the stability of the proposed approach, we compared the results
of simulations with a varying time resolution for the 1D blood flow scheme or the 3D LV model, respectively.
Taking into account the intrinsic time-step limitation necessary to ensure stability of the 1D blood flow scheme, associated with a CFL condition \cite{vos2010h}, we did not consider a time step $\mathrm{dt1D}$ higher than $\SI{1e-3}{\s}$. 
In more detail, in the first test case we considered
$\mathrm{dt1D}\in\{\SI{5e-4}{\s}, \SI{1e-4}{\s}, \SI{5e-5}{\s}\}$,
while keeping the time resolution for the 3D LV model constant at
$\mathrm{dt3D}=\SI{1e-3}{\s}$.  
The circulatory system is represented in this case by one vessel segment of $\SI{200}{\mm}$ and an RCR Windkessel model for terminal boundary conditions.
We emphasise that the coupling is always performed at the time resolution of the 3D EM model.
\Cref{fig:res_dt1d} demonstrates that solutions remain unchanged
when considering a different time discretisation of the circulatory model. 
Therefore, our results indicate that we consider the lowest temporal resolution ensuring stability of the 1D numerical scheme.
\newline
As a second example, we performed two simulations considering a varying time
resolution for the 3D heart model, $\mathrm{dt3D} \in \{\SI{1e-3}{\s}, \SI{5e-4}{\s}\}$,
keeping the same temporal resolution $\mathrm{dt1D} = \SI{1e-4}{\s}$ for the 1D
circulatory model.
\Cref{fig:res_dt3d} shows that, also in this case, solutions remain
unchanged when considering a different time discretisation of the 3D heart model.
Thus, we will consider $\mathrm{dt3D} = \SI{1e-3}{\s} $ in what follows.
We refer to~\cite{augustin2016anatomically} for a more detailed convergence
testing of the cardiac model alone.
\subsection{Prediction of pulse wave propagation effects on cardiac function}
\label{sec:res_application}
To test the ability of the coupled model to predict the effect of pulse wave propagation on cardiac electromechanics
we consider three test cases, ranging from an idealised aortic segment with a constant radius to a pathological condition called aortic coarctation and a network composed of 116 arteries.
\subsubsection{Idealised aortic segment with constant radius}
\label{sec:126mm_constant}
\begin{figure}[ht]
	\centering
	\begin{center}
		\includegraphics[width=15cm,height=4cm]{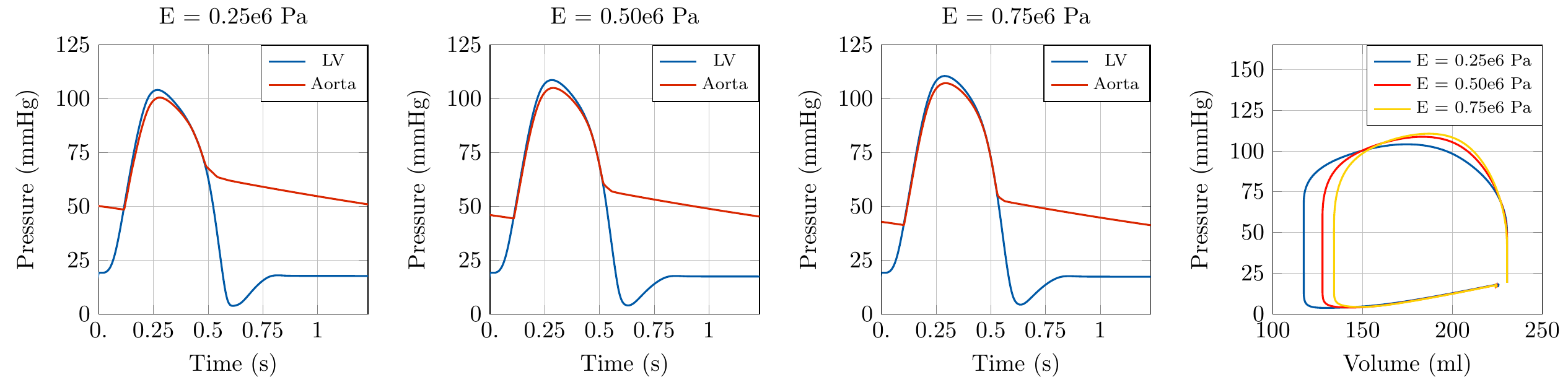}
		\caption{Illustration of model predictions. Test case 1: Idealised aortic segment with constant radius. Left: Pressure trace in the LV and at the inlet of the aorta. Right: Pressure-volume loop in the LV. Comparison of model predictions considering different Young's modulus $E$ in the 1D blood flow model.}
		\label{fig:126_mm_E_p}
	\end{center}
\end{figure}

\begin{figure}[ht]
	\centering
	\begin{center}
		\includegraphics[width=15cm,height=5cm]{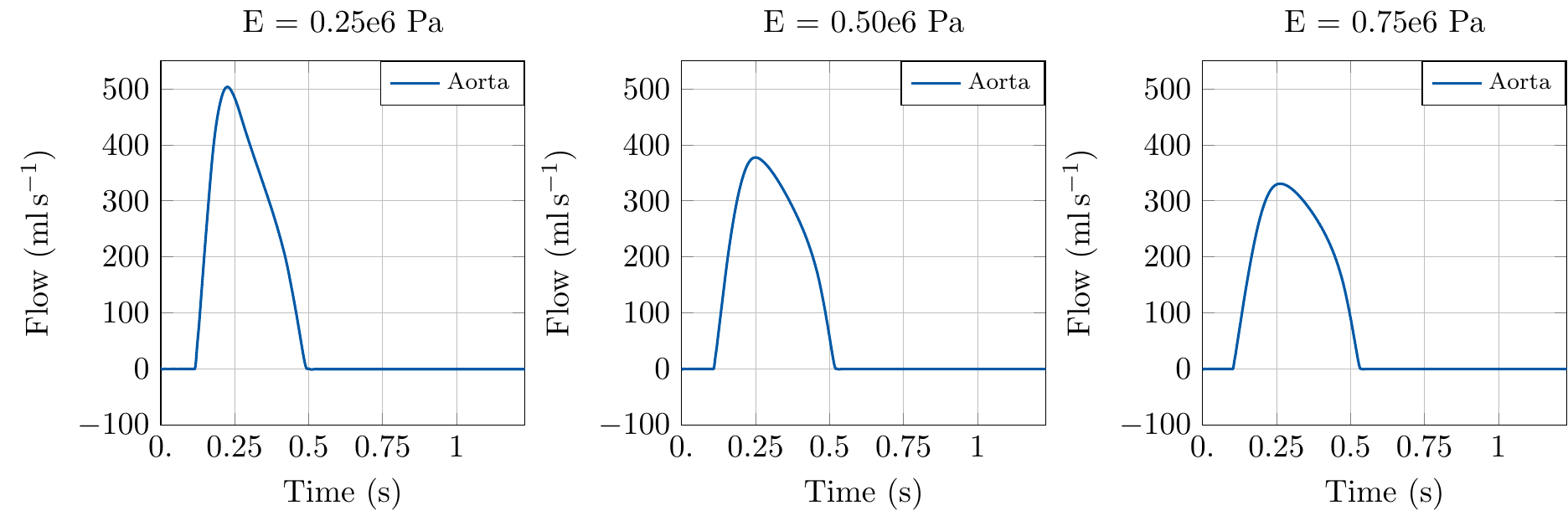}
		\caption{Illustration of model predictions. Test case 1: Idealised aortic segment with constant radius. Flow trace in the LV and at the inlet of the aorta.  Comparison of model predictions considering different Young's modulus $E$ in the 1D blood flow model.}
		\label{fig:126_mm_E_q}
	\end{center}
\end{figure}

The first example illustrates the effect of vessel stiffening on haemodynamics
and LV function. In particular, we considered a 3D LV solid model coupled to
a vessel segment with constant radius and a length of \SI{126}{\mm} with lumped terminal boundary conditions.
For the baseline case we set the vessel wall stiffness described by the Young modulus
to $E = \SI{0.25e6}{\Pa}$.
For two test cases with increased stiffness we set the Young modulus $E$ to
\SI{0.50e6}{\Pa} and \SI{0.75e6}{\Pa}, respectively.
The 1D blood flow model was initialised using an inflow boundary condition. For the simulations in this test case we fixed $\mathrm{dt1D} =\SI{1e-4}{\s}$.
Simulation results of the coupled model are shown in~\Cref{fig:126_mm_E_p}, where the
pressure trace in the LV and at the inlet of the aortic segment is compared for the
different stiffness values.
It can be observed that an increase in $E$ produced an increase
in peak pressure and a substantial change in the shape of the pressure trace.
This is due to the fact that the stiffening of the vessel wall affects pulse wave propagation in the vessel itself.
In particular, the increase of aortic stiffness is associated with a premature rise and decay of the pressure wave and the well-known phenomenon of peak pressure augmentation~\cite{nichols1991mcdonald}. 
The increased stiffness of the aortic vessel has a significant impact on
LV function, see the PV loop in~\Cref{fig:126_mm_E_p} and ~\Cref{fig:126_mm_E_q}.
The stiffening of the arterial vessel is associated with an
increase in end-systolic volume (ESV) at unchanged end-diastolic volume (EDV), leading to a reduction of stroke volume (SV)
\subsubsection{Stenotic aortic segment}
\label{sec:126mm_stenosis}
\begin{figure}[ht]
	\centering
	\begin{center}
		\includegraphics[width=15cm,height=4cm]{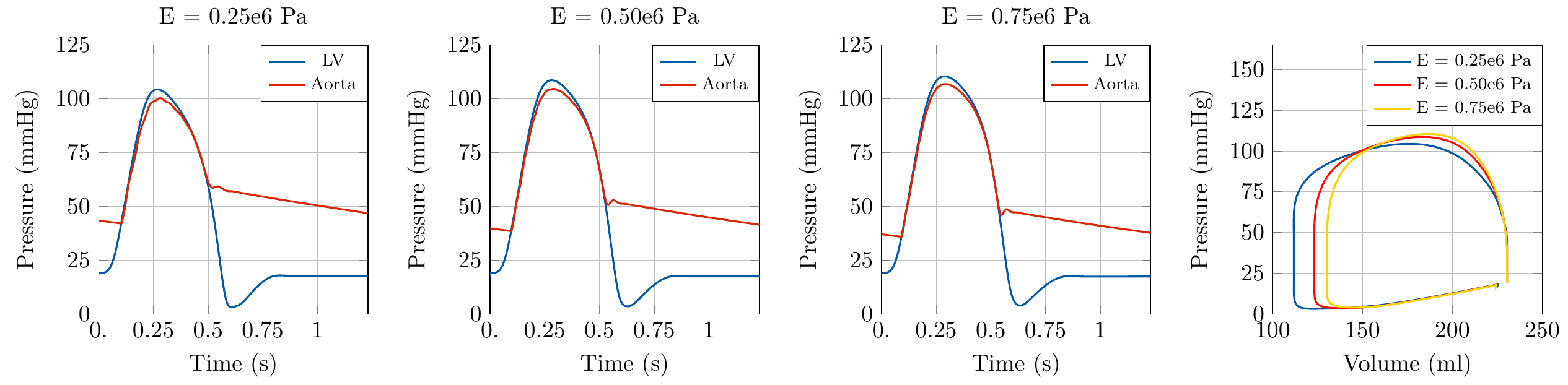}
		\caption{Illustration of model predictions. Test case 2: Aortic segment with coarctation.  Left: Pressure trace in the LV and at the inlet of the aorta. Right: Pressure-volume loop in the LV. Comparison of model predictions considering different Young's modulus $E$ in the 1D blood flow model. }
		\label{fig:CoA_E_p}
	\end{center}
\end{figure}

\begin{figure}[ht]
	\centering
	\begin{center}
		\includegraphics[width=12cm,height=5cm]{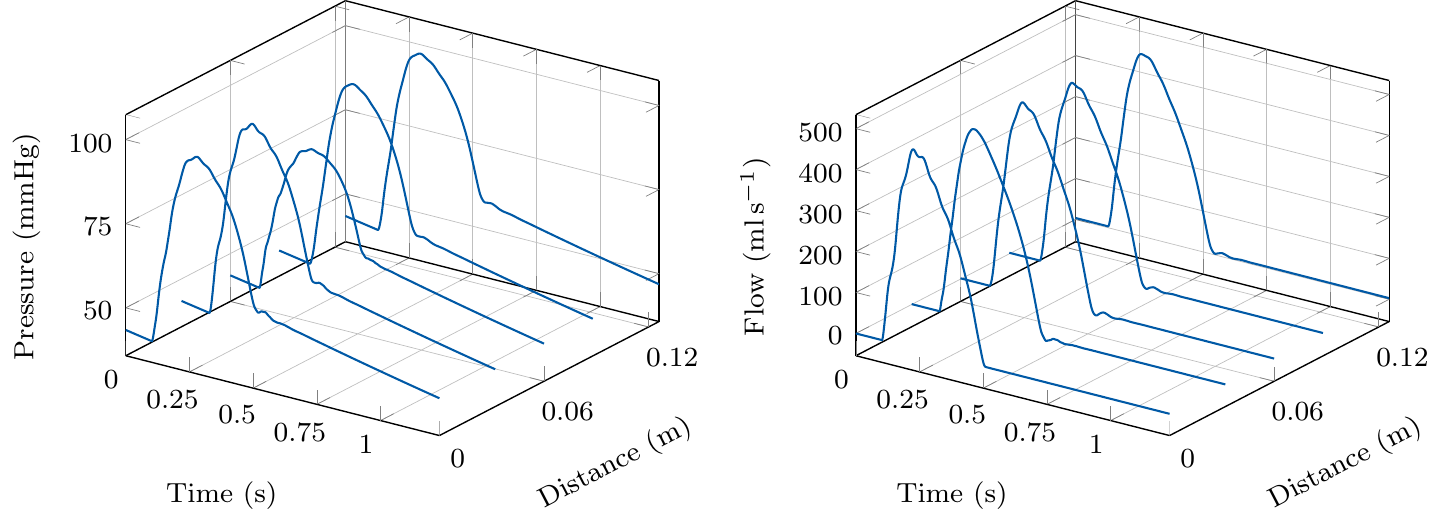}
		\caption{Illustration of model predictions. Test case 2: Aortic segment with coarctation. Pressure and flow profile in the aortic segment. Young's modulus $E = \SI{0.25e6}{\Pa}$. }
		\label{fig:CoA_E_out}
	\end{center}
\end{figure}

As a second example, we considered a pathological configuration of aortic coarctation,
that is a congenital defect consisting of a local narrowing of the aortic lumen
(most commonly in the aortic arch).
To this end, we considered an arterial segment of \SI{126}{\mm} endowed with a stenosis model \cite{jin2021arterial}, with a 30\% stenosis halfway of its length. Terminal boundary conditions are given by an RCR Windkessel model.
In order to explore the effect of stiffening of the vessel in this configuration,
we considered a baseline case with $E = \SI{0.25e6}{\Pa}$
and two test cases with $E$ equal to \SI{0.50e6}{\Pa} and \SI{0.75e6}{\Pa}, respectively.
For the simulations in this test case we set $\mathrm{dt1D} =\SI{1e-5}{\s}$ for stability considerations.
Also in this example, an increase in aortic stiffness is reflected in an
augmentation of peak pressure and a variation in the pressure profile,
see~\Cref{fig:CoA_E_p}.
The effect of vessel stiffness on LV function is shown in the PV loop of~\Cref{fig:CoA_E_p}
consists in a reduction in SV caused by an increase of the ESV.
\begin{figure}[ht]
	\centering
	\begin{center}
		\includegraphics[width=15cm,height=5cm]{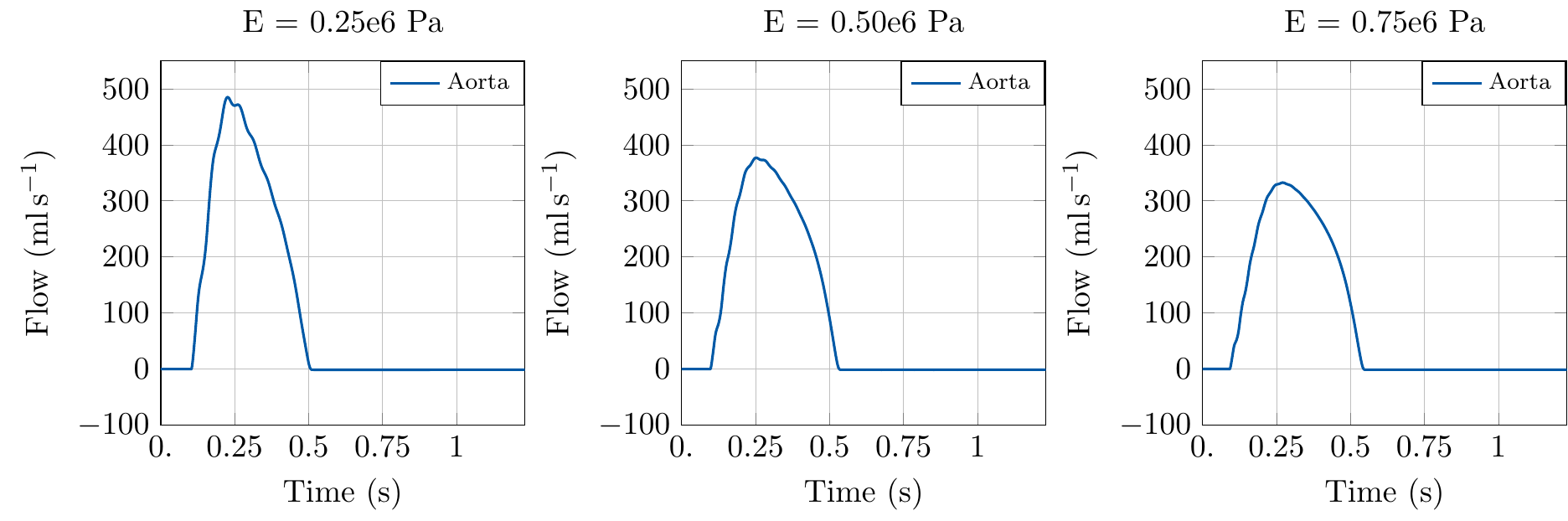}
		\caption{Illustration of model predictions. Test case 2: Aortic segment with coarctation. Flow trace in the LV and at the inlet of the aorta. Comparison of model predictions considering different Young's modulus $E$ in the 1D blood flow model. }
		\label{fig:CoA_E_q}
	\end{center}
\end{figure}
For the sake of completeness we show the variation of pressure and flow signals along the distance of the 1D vessel for the case $E = \SI{0.25e6}{\Pa}$. In particular, we show the model predictions of the 3D-1D coupled model in~\Cref{fig:CoA_E_out}.
\subsubsection{116 aortic segments}
\label{sec:116arteries}
\begin{figure}[ht]
	\centering
	\begin{center}
		\includegraphics[width=15cm,height=5cm]{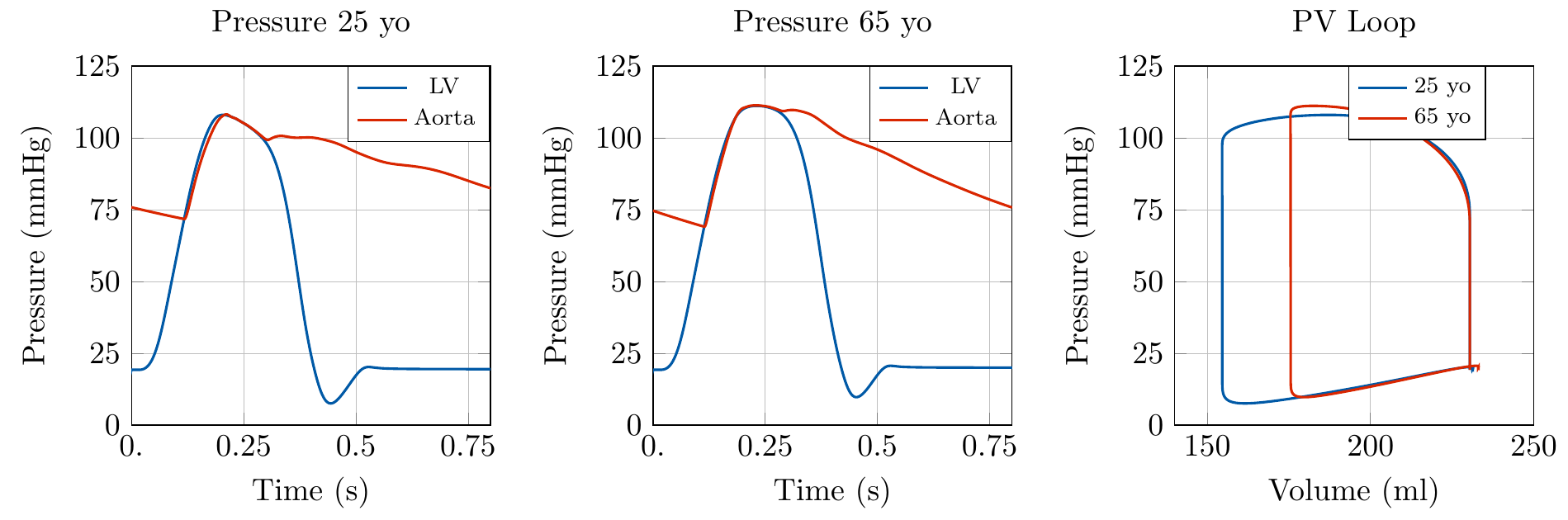}
		\caption{Illustration of model predictions.  Test case 3: Network with 116 vessels. Left: Pressure trace in the LV and at the inlet of the aorta. Right: Pressure-volume loop in the LV. Comparison of model predictions considering different cardiovascular properties in the 1D blood flow model corresponding to a 25 yo and a 65 yo subject, respectively. }
		\label{fig:116_p}
	\end{center}
\end{figure}

\begin{figure}[ht]
	\centering
	\begin{center}
		\includegraphics[width=10cm,height=5cm]{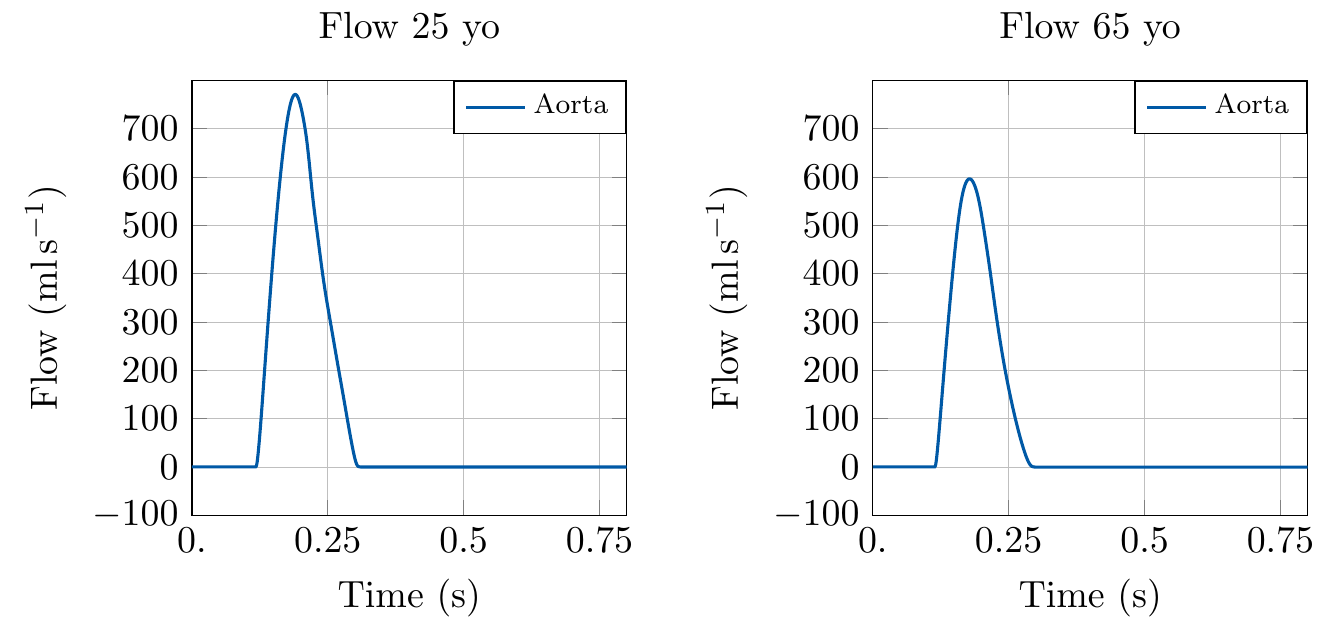}
		\caption{Illustration of model predictions. Test case 3: Network with 116 vessels. Flow trace in the LV and at the inlet of the aorta.  Comparison of model predictions considering different cardiovascular properties in the 1D blood flow model corresponding to a 25 yo and a 65 yo subject, respectively. }
		\label{fig:116_q}
	\end{center}
\end{figure}

\begin{figure}[ht]
	\centering
	\begin{center}
		\includegraphics[width=12cm,height=5cm]{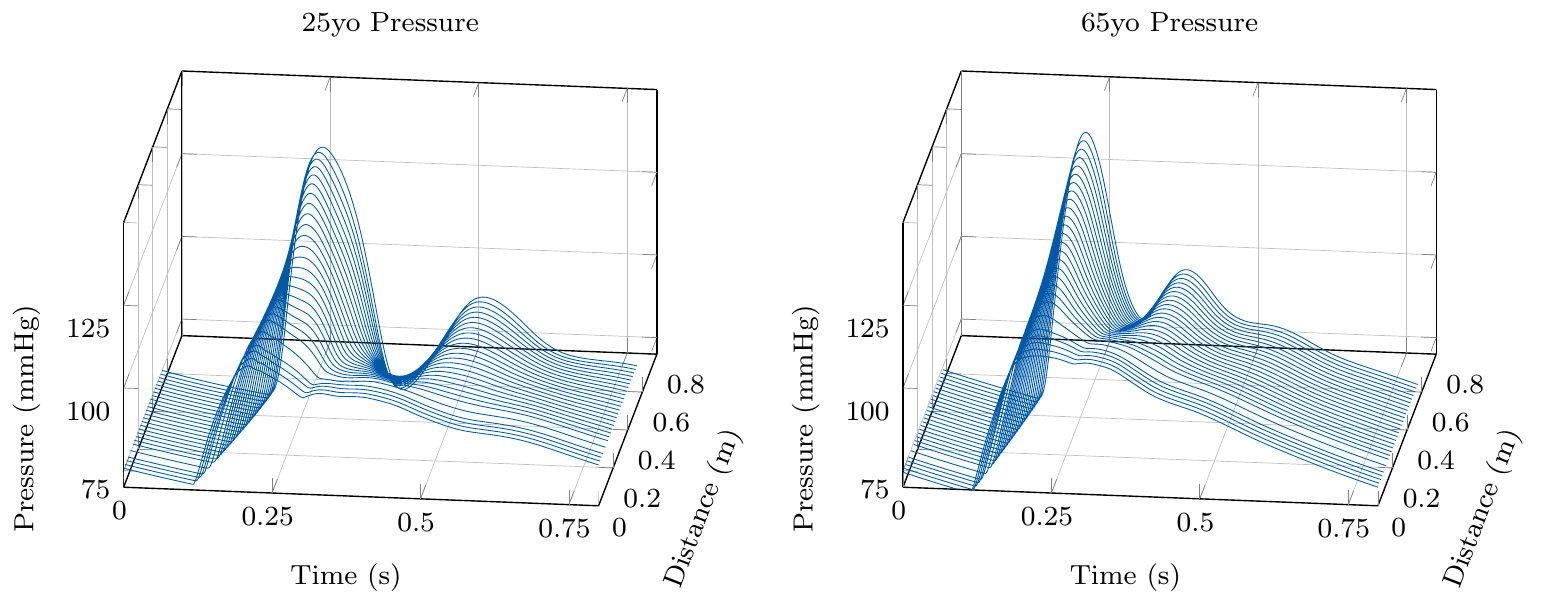}
		\caption{Illustration of model predictions. Test case 3: Network with 116 vessels. Change in pressure trace from the aortic root to the right brachial artery (in the arm). Comparison of model predictions considering different cardiovascular properties in the 1D blood flow model corresponding to a 25 yo and a 65 yo subject, respectively. }
		\label{fig:116_Outp}
	\end{center}
\end{figure}

In order to explore the effect of bifurcations on pulse wave reflections
and their impact on LV function the 3D EM model of the LV was coupled to a more complex model of the arterial system.
To this end, as a third example we coupled the 3D model to a network consisting of 116 arterial vessels. Moreover, we explored two different physiological configurations for the 1D model, corresponding to a 25 year-old subject and a 65 year-old subject. We refer to~\cite{charlton2019modeling} for further detail on the 1D networks considered.  
Concerning the 3D heart model, we considered the same parameter values as in the previous test cases, listed in~\Cref{tab:heart_input_params}, with some exceptions. In particular, the heart cycle period was set to $\SI{800}{\ms} $ in both configurations, and the parameters of the active stress law were adapted according to the new heart cycle period. 
We set $\mathrm{dt1D} =\SI{1e-5}{\s}$ in this test case as well.
Also in this example, the increase in aortic stiffness and systemic vascular resistance associated with healthy ageing is reflected in an
augmentation of the peak pressure and a variation in the pressure profile,
see~\Cref{fig:116_p}.
The effect of ageing on LV function is also striking, as shown in ~\Cref{fig:116_p,fig:116_q}.
In particular, we can see a reduction in SV caused by an increase of the ESV.
For this test case we can show the effect of pulse pressure amplification, i.e. the amplification in pulse pressure signal towards more distal locations, associated with reflections at bifurcation sites. For the sake of illustration, in~\Cref{fig:116_Outp} we consider the change in pressure trace from the aortic root to the right brachial artery. As expected, pulse pressure amplification was more pronounced in the younger subject than in the 65-yo subject.
\section{Discussion}%
\label{sec:discussion}
In this methodological study we have developed a numerical scheme
that allows to couple an advanced 3D EM model of cardiac function to a 1D model of the human circulation.
First, we have demonstrated the feasibility of the proposed approach by coupling a
patient-specific 3D LV model to a 1D model of the aorta and shown the ability of the
model to correctly predict known physiological behaviours associated with
arterial wall stiffening both in a healthy case and in the pathological case of aortic coarctation.
Furthermore, we have shown an application to the coupling of the 3D LV model with a complex network composed by 116 arteries, in order to study the effect of bifurcations on pulse wave reflections and their impact on cardiovascular function.
\subsection{Prediction of pulse wave propagation effects on cardiac function}
The results in \Cref{sec:results} show an augmentation of peak pressure in the LV and pulse pressure in the aorta and a
substantial variation in the shape of the pressure trace as a consequence of
stiffer vessel walls.
Peak pressure in the LV increases with increasing E, since there is a higher afterload.
In the aorta, the pulse pressure increases with increasing Young's modulus E, since the aortic model becomes less compliant.
Physiologically, the pulse pressure trace is composed of two additive components: first,
a forward propagating wave generated by LV contraction and second,
a backward propagating wave from the periphery.
From this, the behaviour observed in the simulations is to be expected as the
increased aortic stiffness generates a faster pulse wave propagation.
Hence, this causes a premature arrival of the forward propagating wave and, in turn, of the backward propagating wave,
entailing peak pressure augmentation~\cite{nichols1991mcdonald}.\\
In addition, the stiffening of the arterial vessel is related to
an increase in ESV.
The main reason for this effect is that a lower compliance
of the vessel, i.e., its ability to distend and increase volume with increasing
transmural pressure, is responsible for a reduction of blood ejected by the
heart chamber, thus producing a greater ESV.
Since EDV is influenced by the cardiac passive behaviour and preload rather than
by afterload, this volume does not vary when the aortic vessel properties change.
Hence, the stiffening of the aortic vessel leads to a reduction in SV.
We have shown that the proposed 3D-1D model is capable to reproduce the impact
of increased stiffness and pulsatile load in the circulation and their effect
on heart performance.
These results confirm that such models are particularly suited to explore the effects
of changes in pulse wave propagation, which can hardly be captured using
lumped parameter models~\cite{Alastruey2012a}.
Also, effects of these interactions are very complex to study \emph{in vivo}
due to the need for high fidelity and simultaneous measurements of cardiac function
and circulation at several locations, together with the capability to assess chamber
interactions~\cite{Mynard_Smolich_2015}.
Hence, 3D-1D models are of high relevance to investigate
interactions between vascular waves and cardiac chamber function.
Overall, a coupling to 1D models should be preferred over a coupling to
lumped parameter models for clinical applications that are profoundly influenced
by disrupted blood flow propagation, e.g., aortic coarctation pathologies and
pulmonary arterial hypertension.
\subsection{Numerical aspects}
Computationally, the 3D-1D approach is comparable to 3D-0D models, as the extra cost of solving the 1D model is essentially negligible compared to the cost of the 3D model.
As such, whenever wave transmission effects are to be investigated, the 3D-1D model is preferable as it can be used instead of a 3D-0D model as in \cite{augustin2021computationally} without any computational penalties.
Further, we emphasise that the implementation of the coupling framework is
not constrained to the specific choice of the 1D solver used for the solution
of the blood flow equations in the circulatory system.
In personalised applications the extra cost of identifying parameters of a distributed 1D system must be factored in though.
\subsection{Limitations}
An extension of the study is the use of more complex biventricular and
four-chamber EM cardiac models \cite{strocchi2020publicly}, in order to consider the interactions among the
heart chambers and thus provide more physiological simulations of the cardiac function.
Also, it is well known that cardiac preload is mainly affected by venous return.
As a consequence, in order to study more accurately the feedback of circulation (including venous return) on the cardiac function, the use of a 1D global model of the human circulation is more appropriate.
Thus, another future perspective consists in coupling the 3D EM cardiac model with closed-loop global 1D blood flow models, e.g.~\cite{muller2014global}, in order to explore the effect of circulation on preload as well.
We emphasise that the coupling strategy holds for a general framework consisting of multiple cardiac cavities and is prone to different discretisations and geometries for the circulatory model, therefore such extensions are feasible.
Third, this study did not specifically focus on the parameterisation of the
coupled model to patient-specific data. However, the personalisation of
cardiovascular models is of crucial importance for their predictive role,
since the circulatory system has a complex network structure and shows a
significant inter-individual variability.
Owing to their distributed nature, 1D models of blood flow involve a larger
amount of parameters that must be calibrated than lumped parameter models, related to geometric aspects,
elastic properties and peripheral impedances of the
network~\cite{Reymond_etal_2011,Blanco_etal_2014b}.
This poses serious limitations on the applicability of the use of detailed 1D
circulatory models in clinical applications.
Nonetheless, there are numerous contributions in the scientific literature to cope
with this problem, in terms of simplification of the geometry and reduction of the
parameter set, together with the application of optimisation methods and robust
inverse problem strategies.
For example, the topological complexity of the arterial tree can be optimised by
effectively reducing the number of arterial segments included in the 1D model,
still preserving the key characteristics of flow and pressure
waveforms~\cite{epstein2015reducing,fossan2018optimization}.
Common approaches also contemplate the rescaling of distal properties and vessel
geometry based on allometric scales and global descriptors related to the geometry and properties of the
network that are easy to measure, such as pulse wave velocity, body size,
and biological age~\cite{quarteroni2017integrated}.
The level of detail of the network representing the circulatory system should be carefully considered taking into account the specific clinical application under study. \\
In addition, recent parameter estimation methods based on data assimilation techniques,
i.e., algorithms combining mathematical models with available measurements to improve
the accuracy of model predictions and estimate patient-specific parameters,
have provided promising results in cardiovascular
applications~\cite{lombardi-14,caiazzo2017assessment,muller2018reduced}.
\section{Conclusion}%
\label{sec:conclusion}
We developed a stable strategy to perform a coupling between
a 3D EM model of the heart and a 1D blood flow model of the arterial system, based on the resolution of a saddle-point problem for the volume
and pressure in the cavity.
We built a personalised coupled 3D-1D model of LV and arterial system and showed robustness and accuracy of our scheme in a numerical benchmark. 
We demonstrated the ability of the coupled system to efficiently simulate physiological response to alterations in features of the circulatory system influencing pulse wave transmission, including aortic stiffening, aortic stenosis and bifurcations causing wave reflections.
The additional computational costs associated with the use of 1D models instead of standard 0D models in the coupled system are negligible. 
As a consequence, the use of our coupled 3D-1D model is beneficial for a broad spectrum of clinical applications where wave transmission effects are under study,
such as aortic coarctation and pulmonary arterial hypertension.   
\section*{Acknowledgements}
The study received support from the Austrian Science Fund (FWF) grant I2760-B30 and from BioTechMed Graz grant ILEARNHEART\@.
Simulations for this study were performed on the Vienna Scientific Cluster (VSC-4), which is maintained by the VSC Research Center in collaboration with the Information Technology Solutions of TU Wien. The authors acknowledge the financial support by the University of Graz. The authors acknowledge Dr. Laura Marx, PhD (Medical University of Graz, Austria) for the technical support with the calibration of the 3D EM cardiac model.
\bibliographystyle{elsarticle-num-names}
\bibliography{../MyBibFile}
\clearpage
\begin{appendix}
\section{Electromechanical PDE model}%
\label{sec:app_heart}
	This section is devoted to the description of the mathematical models considered to describe the most fundamental aspects of cardiovascular function, comprising electrophysiology, passive and active mechanics and haemodynamics.
	\subsection{Electrophysiology}
	A recently developed R-E model~\cite{neic17:_efficient} is
	considered for the generation of electrical activation sequences
	serving as a trigger for active stress generation in cardiac tissue.
	This hybrid R-E model combines a standard R-D model
	based on the monodomain equation with an eikonal model.
	In more detail, we consider the eikonal equation given as
	\begin{equation} \label{eq:_eikonal}
	  \left\{
	    \begin{array}{rcll}
	      \sqrt{\nabla_{\vX} t_\mathrm{a}^\top \, \tensor{V} \,
	      \nabla_{\vX} t_\mathrm{a}} & =
	      & 1 \qquad & \text{in }  \Omega_0, \\
	      t_\mathrm{a} & = & t_0 & \text{on } \Gamma_0^{\ast},
	    \end{array}
	  \right.
	\end{equation}
	where $(\nabla_{\vX})$ is the gradient with respect to the end-diastolic
	reference configuration $\Omega_{0}$,
	$t_\mathrm{a}$ is a positive function that describes the wavefront arrival
	time at location $\vX \in \Omega_0$,
	and $t_0$ denote the initial activations at locations
	$\Gamma_0^\ast \subseteq \Gamma_{0,\mathrm{N}}$.
	The symmetric positive definite $3 \times 3$ tensor $\tensor{V}(\vX)$
	contains the squared velocities
	$\left(v_{\vf}(\vX),v_{\vs}(\vX),v_{\vn}(\vX)\right)$
	associated with the tissue's eigenaxes $\vf_0$, $\vs_0$, and $\vn_0$.
	Then, the arrival time function $t_\mathrm{a}(\vX)$ is used in a
	modified monodomain R-D model given as
	\begin{equation} \label{equ:_monodomain_R_E}
	  \beta C_\mathrm{m} \frac{\partial V_\mathrm{m}}{\partial t} =
	    \nabla_{\vX} \cdot \boldsymbol{\sigma}_\mathrm{i} \nabla_{\vX} V_\mathrm{m} +
	    I_\mathrm{foot} - \beta I_\mathrm{ion},
	\end{equation}
	where an arrival time dependent foot current,
	$I_{\mathrm{foot}}(t_\mathrm{a})$, is added,
	in order to mimic sub-threshold electrotonic currents
	and produce a physiological foot of the action potential.
	The key advantage of this R-E model is that it enables
	to compute activation sequences at much coarser spatial resolutions
	without being afflicted by the spatial undersampling artefacts
	leading to conduction slowing or even numerical conduction blocks that are
	observed in standard R-D models.
	The tenTusscher--Noble--Noble--Panfilov model %
	of the human ventricular myocyte~\cite{tentusscher04:_TNNP} is employed to model Ventricular EP.
	\begin{table}[htbp]
		\centering
		\footnotesize
		\begin{tabularx}{0.85\textwidth}{lrlL}
			\toprule
			Parameter & Value & Unit & Description \\
			\midrule
			\multicolumn{4}{l}{\emph{Passive biomechanics}}\\
			$\rho_0$        & \num{1060.0}  & \si{\kg/\m\cubed} & tissue density\\
			$\kappa$        & 650  & \si{\kPa} & bulk modulus\\
			$a$             & 0.8  & \si{\kPa} & stiffness scaling\\
			$b_\mathrm{ff}$ & 5.0  & [-] & fibre strain scaling\\
			$b_\mathrm{ss}$ & 6.0  & [-] & cross-fibre in-plain strain scaling\\
			$b_\mathrm{nn}$ & 3.0  & [-] & radial strain scaling \\
			$b_\mathrm{fs}$ & 10.0 & [-] & shear strain in fibre-sheet plane scaling \\
			$b_\mathrm{fn}$ & 2.0  & [-] & shear strain in fibre-radial plane scaling \\
			$b_\mathrm{ns}$ & 2.0  & [-] & shear strain in transverse plane scaling \\
			\midrule
			\multicolumn{4}{l}{\emph{Active biomechanics}}\\
			$\lambda_0$           & 0.7   & \si{\ms}    &  onset of contraction\\
			$V_\mathrm{m,Thresh}$ & -60.0 & \si{\milli\volt} & membrane potential threshold \\
			$t_\mathrm{emd}$      & 15.0  & \si{\ms}    &  EM delay\\
			$S_\mathrm{peak}$     & 60    & \si{\kPa} &  peak isometric tension\\
			$t_\mathrm{dur}$      & 575.0 & \si{\ms}    &  duration of active contraction\\
			$\tau_{c_0}$          & 105.0 & \si{\ms}    &  baseline time constant of contraction\\
			$\tau_\mathrm{r}$     & 90.0 & \si{\ms}    &  time constant of relaxation\\
			$\mathrm{ld}$         &   35.0 & [-]         &  degree of length-dependence\\
			$\mathrm{ld}_\mathrm{up}$ & 100.0 & \si{\ms}&  length-dependence of upstroke time\\
			\midrule
			\multicolumn{4}{l}{\emph{Electrophysiology}}\\
			$t_\mathrm{cycle}$ & \num{1.231}  & \si{\s}
			& cycle time ($=1/\mathrm{heartrate}$)\\
			AA delay & 20.0   & \si{\ms} & inter-atrial conduction delay \\
			AV delay & 100.0  & \si{\ms} & atrioventricular conduction delay \\
			VV delay & 0.0    & \si{\ms} & inter-ventricular conduction delay \\
			$(v_\mathrm{f},v_\mathrm{s},v_\mathrm{n})$ & (0.6, 0.4, 0.2)    & \si{\m/\s} & conduction velocities \\
			$(g_\mathrm{f},g_\mathrm{s},g_\mathrm{n})$ & (0.44, 0.54, 0.54)  & \si{\m/\s} & conductivities in LV \\ %
			$\beta$ & 1/1400 & \si{\cm^{-1}} & membrane surface-to-volume ratio \\
			$C_\mathrm{m}$ & 1 & \si{\micro\farad/cm^2} & membrane capacitance \\
			\bottomrule
		\end{tabularx}
		\caption{Input parameters for the 3D PDE model of the left ventricle (\LV).
			Adjusted to match patient-specific data.}%
		\label{tab:heart_input_params}
	\end{table}
	\subsection{Passive mechanics of the heart}
	\label{sec:app_passive_heart}
\begin{figure}[ht!]
	  \centering
	  \includegraphics[width=1.0\linewidth,keepaspectratio]{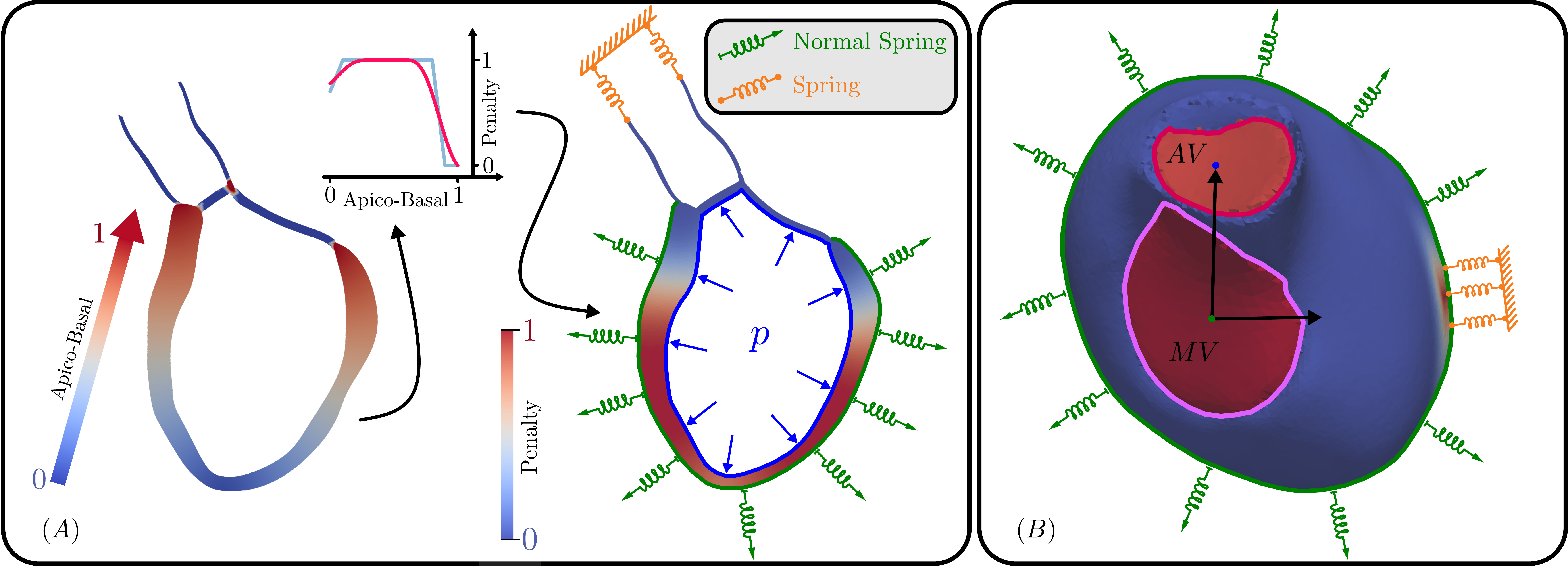}
	  \caption{Boundary conditions applied to the LV models.}%
	  \label{fig:mechbc}
	\end{figure}

	The deformation of the heart is governed by imposed external loads such as pressure in the cavities or from surrounding tissue and active stresses intrinsically generated during contraction.
	The cardiac tissue is characterised as a hyperelastic, nearly incompressible, anisotropic material with a nonlinear
	stress-strain relationship. Mechanical deformation is described by Cauchy's equation of motion
	leading to a  boundary value problem
	\begin{equation}\label{eq:bvp}
	  \rho_0\,\ddot{ \vu}(t,\vX)-\nabla_{\vX}\cdot\mathbf{F}\mathbf{S}(\vu,\vX) = \vb_0(\vX)
	  \quad \mbox{in }  \Omega_0 ,
	\end{equation}
	for $t \in [0, T]$, where $\rho_0$ is the density in the Lagrange reference configuration,
    $\vu(t,\vX)$ is the unknown nodal displacement,  $\dot{\vu}(t,\vX)$ is the nodal velocity,  $\ddot{\vu}(t,\vX)$ is the nodal acceleration,
	$\mathbf{F}(\vu,\vX)$ is the deformation gradient,
	$\mathbf{S}(\vu, \vX)$ is the second Piola--Kirchhoff stress tensor, $\vb_0(\vX)$ represents the body forces  %
	and $(\nabla_{\vX}\;\cdot)$ denotes the divergence operator in the reference configuration.
	We consider as initial conditions
		\begin{equation*}
		\vu(0,\vX) = \vZero, \quad  \dot{\vu}(0,\vX) = \vZero,
		\end{equation*}
		and we set $\vb_0(\vX) = \vZero$ for the sake of simplicity.
	The boundary of the left ventricular model is decomposed in three components, i.e., $\partial\Omega_0=\overline{\Gamma}_{\mathrm{endo},0} \cup
	\overline{\Gamma}_{\mathrm{epi},0}\cup \overline{\Gamma}_{\mathrm{base},0}$,
	where $\overline{\Gamma}_{\mathrm{endo},0}$ denotes the endocardium, $\overline{\Gamma}_{\mathrm{epi},0}$ the epicardium, and
	$\overline{\Gamma}_{\mathrm{base},0}$ the base of the ventricle.
At the endocardium normal stress boundary conditions are imposed:
\begin{equation}\label{eq:neumann}
\tensor{F}\tensor{S}(\vec{u},\vec{X})\,\normalout(\vec{X})
= -p(t) J\tensor{F}^{-\top}\normalout(\vec{X})
\quad\text{on}\quad {\Gamma}_{\mathrm{endo},0}\times (0,T)
\end{equation}
with $\normalout$ the outer normal vector.
Omni-directional spring type boundary conditions
	constrain the ventricle at the basal cut plane $\overline{\Gamma}_{\mathrm{base},0}$~\cite{land2018influence},
	and spatially varying normal Robin boundary conditions
	are applied at the epicardium $\overline{\Gamma}_{\mathrm{epi},0}$~\cite{strocchi2020simulating} to simulate the pericardial mechanical constrains.\\
	We consider the following decomposition of the total stress $\mathbf{S}$:
	\begin{equation} \label{eq:additiveSplit}
	  \mathbf{S}= \mathbf{S}_\mathrm{pas}+ \mathbf{S}_\mathrm{act},
	\end{equation}
	where $\mathbf{S}_\mathrm{pas}$ and $\mathbf{S}_\mathrm{act}$
	denote the passive and active stresses, respectively.
	Passive stresses are modelled according to the constitutive law
	\begin{equation}
	  \mathbf{S}_\mathrm{pas}=2\frac{\partial\Psi(\mathbf{C})}{\partial\mathbf{C}}
	\end{equation}
	with $\Psi$ an invariant-based strain-energy function used to model the anisotropic behaviour of cardiac tissue.
	Numerous constitutive models have been developed in the literature, ranging from simpler transversely-isotropic models to more complex orthotropic laws.
	In this article, we consider a hyper-elastic and transversely isotropic strain-energy function $\Psi$
	\begin{equation} \label{eq:guccioneStrainEnergy}
            \Psi_{\mathrm{Guc}}(\mathbf{C}) = \frac{\kappa}{2} {\left( \log\,J \right)}^2 +
	  \frac{C_\mathrm{Guc}}{2}\left[\exp(\mathcal{Q})-1\right]
	\end{equation}
	proposed by Guccione et al.~\cite{guccione1995finite}.
	Here, the term in the exponent is
	\begin{equation}
	  \label{eq:guccioneQ}
	  \begin{split}
	  \mathcal{Q} = \; &
            b_{\mathrm{f}} {(\vf_0\cdot\overline{\mathbf{E}}\,\vf_0)}^2 +
            b_{\mathrm{t}} \left[{(\vs_0\cdot\overline{\mathbf{E}}\,\vs_0)}^2+
                                 {(\vn_0\cdot\overline{\mathbf{E}}\,\vn_0)}^2+
                             2{(\vs_0\cdot\overline{\mathbf{E}}\,\vn_0)}^2\right]+ \\
             & 2b_{\mathrm{fs}} \left[{(\vf_0\cdot\overline{\mathbf{E}}\,\vs_0)}^2+
             {(\vf_0\cdot\overline{\mathbf{E}}\,\vn_0)}^2\right]
	  \end{split}
	\end{equation}
	and $\overline{\mathbf{E}}=\frac{1}{2}(\overline{\mathbf C}-\mathbf{I})$
	is the modified isochoric Green--Lagrange strain tensor,  %
	where $\overline{\mathbf C} := J^{-2/3} \mathbf{C}$ with $J=\det\mathbf{F}$.
	Default values of $b_{\mathrm{f}}=18.48$, $b_\mathrm{t}=3.58$,
	and $b_\mathrm{fs}=1.627$ are used.
	The parameter $C_\mathrm{Guc}$ is used to fit the LV model to
	an empirical Klotz relation~\cite{klotz2007computational} by a combined
	unloading and re-inflation procedure.
	The bulk modulus $\kappa$, which serves as a penalty parameter
	to enforce nearly incompressible material behaviour, is set as
	$\kappa = \SI{650}{\kPa}$.
	\paragraph{Mechanics boundary conditions at the LV}
	The different BCs applied to the LV models are summarised in Fig.~\ref{fig:mechbc}.
	The springs attached to the aortic rim and at the pericardium are shown
	in Fig.~\ref{fig:mechbc}A, as well as the pressure BC in the cavity at the
	endocardial surface. The pericardial springs penalise displacement only in normal direction and are gradually scaled from the apex to the base. Therefore,
	the distance in apico-basal direction is used to create a penalty map, see Fig.~\ref{fig:mechbc}A.
	To avoid non-physiological rotation, further springs are attached to the septum,
	see Fig.~\ref{fig:mechbc}B. The location of the septal springs is selected automatically by constructing a local coordinate system spanned by the centres of the apical region, the MV and the AV\@.
	\subsection{Active mechanics}
	\label{sec:app_active_mechanics}
	Following the approach proposed by~\cite{Genet2014, walker2005mri},  we assume that stresses due to active contraction are orthotropic, with full contractile force along the myocyte fibre orientation $\vf_0$ and
	\SI{40}{\%} contractile force along the sheet orientation $\vs_0$.
	In more detail, we define the active stress tensor $\tensor{S}_{\mathrm{a}}$ as
	\begin{equation} \label{eq:act}
	  \tensor{S}_{\mathrm{a}}
            = S_{\mathrm{a}} {\left(\vf_0\cdot\tensor{C}\vf_0\right)}^{-1}
	      \vf_0 \otimes \vf_0
            + \num{0.4}\,S_{\mathrm{a}} {\left(\vs_0\cdot\tensor{C}\vs_0\right)}^{-1}
	      \vs_0 \otimes \vs_0,
	\end{equation}
	with $S_\mathrm{a}$ the scalar active stress that describes the contractile
	force.
	Active stress generation is modelled using a simplified phenomenological contractile model~\cite{niederer11:_length}.
	Owing to its small number of parameters and its
	direct relation to clinically measurable quantities such as peak pressure
	and the maximum rate of rise of pressure, this model is easy to fit
	and thus suitable for being used in clinical EM modelling studies.
	In particular, the active stress transient is defined as
	\begin{equation} \label{eq:_tanh_stress}
	  S_\mathrm{a}(t,\lambda) = S_\mathrm{peak} \,
	  \phi(\lambda) \,
	  \tanh^2 \left( \frac{t_{s}}{\tau_\mathrm{c}} \right) \,
	  \tanh^2 \left( \frac{t_\mathrm{dur} - t_{\mathrm{s}}}{\tau_\mathrm{r}} \right),
	  \qquad \text{for } 0 < t_{\mathrm{s}} < t_\mathrm{dur},
	\end{equation}
	where
	\begin{equation} \label{equ:_tanh_stress2}
	  \phi = \tanh (\mathrm{ld} (\lambda - \lambda_0)),\quad
	  \tau_\mathrm{c} = \tau_\mathrm{c_0} + \mathrm{ld}_\mathrm{up}(1-\phi),\quad
	  t_{\mathrm{s}} = t - t_\mathrm{a} - t_\mathrm{emd}
	\end{equation}
	and $t_{\mathrm{s}}$ is the onset of contraction,
	$\phi (\lambda)$ is a nonlinear length-dependent function
	in which $\lambda$ is the fibre stretch and
	$\lambda_0$ is the lower limit of fibre stretch below which no further active
	tension is generated,
	$t_{\mathrm{a}}$ is the local activation time from Eq.~\eqref{eq:_eikonal},
	$t_{\mathrm{emd}}$ is the EM delay between the onsets of
	electrical depolarisation and active stress generation,
	$S_{\mathrm{peak}}$ is the peak isometric tension,
	$t_\mathrm{dur}$ is the duration of active stress transient,
	$\tau_{\mathrm{c}}$ is time constant of contraction,
	$\tau_{\mathrm{c_0}}$ is the baseline time constant of contraction,
	$\mathrm{ld}_{\mathrm{up}}$ is the length-dependence of $\tau_{\mathrm{c}}$,
	$\tau_{\mathrm{r}}$ is the time constant of relaxation,
	and $\mathrm{ld}$ is the degree of length dependence.
	Thus, active stresses are only length-dependent in this simplified model,
	but dependence on fibre velocity $\dot{\lambda}$ is neglected.
\section{Circulatory system: 1D PDE model}%
\label{sec:oneDmodelling}
	In what follows, the main mathematical aspects of 1D blood flow modelling in the human circulation are detailed.
	\subsection{Some features of 1D blood flow modelling}
	The 1D equations of blood flow can be derived by integrating the Navier-Stokes equations of a generic section of the tube (after assuming axisymmetric flow in a cylindrical tube of radius $R$), obtaining
	\begin{equation}
	\begin{cases}
	\displaystyle
	\frac{\partial A}{\partial t}+\frac{\partial Q}{\partial x}=0 \\[6pt]
	\displaystyle
	\frac{\partial Q}{\partial t} + \frac{\partial }{\partial x}\Big(\alpha \frac{Q^2}{A} \Big)+\frac{A}{\rho} \frac{\partial P}{\partial x}= \frac{f}{\rho},
	\end{cases}
	\label{eq:1DAlastruey}
	\end{equation}
	where $A(x,t)$ is the \textit{cross-sectional area} of the lumen, $Q(x,t)$ is the \textit{average flux} and $P(x,t)$ is the \textit{average internal pressure} over the cross-section.
	Note that in the momentum balance a momentum flux correction factor $\alpha$ (also known as Coriolis coefficient) is introduced, defined as
	\begin{equation*}
	\alpha = \frac{1}{A\, \bar u^2}\int_{S(x,t)} u_x^2(x,t) \, \der \sigma
	\end{equation*}
	An energy inequality that bounds a measure of the energy of the hyperbolic system was derived in~\cite{Formaggia2001}.
	Alternatively, it is possible to rewrite the system in terms of the variables $(A,u)$. In more detail, if we assume a flat velocity profile, i.e. $\alpha = 1$, we obtain the system
	\begin{equation}
	\begin{cases}
	\displaystyle
	\frac{\partial A}{\partial t}+\frac{\partial (Au)}{\partial x}=0 \\[6pt]
	\displaystyle
	\frac{\partial u}{\partial t} + u\frac{\partial u}{\partial x}+\frac{1}{\rho} \frac{\partial P}{\partial x}= \frac{f}{\rho A},
	\end{cases}
	\label{eq:1DAlastrueyAu}
	\end{equation}
	where $A(x,t)$ is the \textit{cross-sectional area} of the lumen, $u(x,t)$ is the \textit{average axial velocity} and $P(x,t)$ is the \textit{average internal pressure} over the cross-section.
	\begin{figure}[ht]
		\centering
		\includegraphics[height=3.cm,width=6.cm]{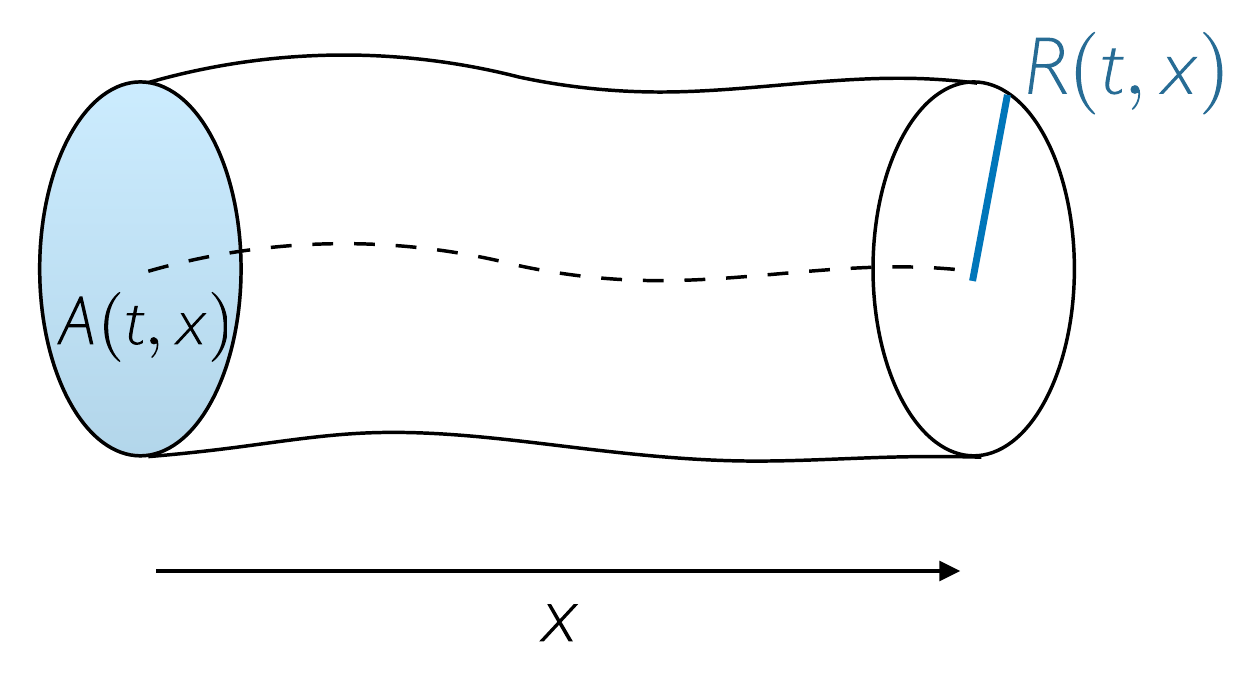}
		\caption{Description of a 1D compliant arterial segment with properties described by a single axial coordinate x.}
	\end{figure}
	For a constant velocity profile satisfying the no-slip condition, the friction force per unit length is
	$f(x,t)= -\eta \,u(x,t)$,
	where $\eta$ is a coefficient depending on the blood viscosity and the \textit{average axial velocity}.
	As System~\eqref{eq:1DAlastrueyAu} is composed by two equations and three unknowns ($A,u,P$), a closure condition is needed, namely the so-called \textit{tube law}:
	\vspace{-0.1cm}
	\begin{equation}
	P(x,t)=P_{ext}(x,t)+K\,\phi\bigl(A(x,t),A_0(x)\bigr)+G \,\Psi\bigl(A(x,t)\bigr),
	\label{eq:TubeLaw}
	\end{equation}
	where $K$ and $G$ are parameters depending on the \textit{wall stiffness} and the wall \textit{viscosity}, respectively.\\
	A simple, purely elastic tube law consists in defining
	\begin{equation*}
	\phi = \frac{\sqrt{A}- \sqrt{A_0}}{A_0}, \ K = \frac{4}{3}\sqrt{\pi}Eh,\  \text{and}\ \Psi =0,
	\end{equation*}
	with $E$ the Young modulus and $h$ the wall thickness of the vessel.
        A more complex tube law is obtained when visco-elastic effects are considered, e.g., by setting
	\begin{equation*}
	\Psi = \frac{1}{A_0\sqrt{A}} \frac{\partial A}{\partial t}, \ G = \frac{2}{3} \sqrt{\pi}\varphi h,
	\end{equation*}
	with $\varphi(x)$ the wall viscosity. In what follows we consider a visco-elastic tube law, see table~\ref{tab:bf_input_params} for more detail.\\
	In the first test cases proposed in this work we consider an arterial network consisting of a segment of the human upper thoracic aorta,  equipped with a 3-element Windkessel model as a terminal boundary condition. A more complex network of vessels is also considered  in~\Cref{sec:116arteries}, and it is taken from~\cite{charlton2019modeling}.\\
		\begin{table}[htbp]
			\centering
			\footnotesize
			\begin{tabularx}{0.7\textwidth}{lrlL}
				\toprule
				Parameter & Value & Unit & Description \\
				\midrule
				$\rho_0$        & \num{1060.0}  & \si{\kg/\m\cubed} & fluid density\\
				$\nu$        & 4e-3  & \si{\Pa\s} & blood flow viscosity\\
				$\alpha$             & 1.1  & [-] & Coriolis coefficient (momentum equation)\\
				$G$ & 1  & \si{\kPa} & vessel wall viscosity\\
				\bottomrule
			\end{tabularx}
			\caption{Input parameters for the 1D PDE model of the arterial circulation.
				}%
			\label{tab:bf_input_params}
		\end{table}
	\subsubsection{Boundary conditions}
	Standard inlet boundary conditions for 1D blood flow models consist of lumped 0D models of the heart or, alternatively, pressure or flow profiles~\cite{Formaggia_etal_2003, Alastruey2012a,muller2014global,Blanco_etal_2014b}.
	Every arterial 1D model must be truncated after a certain number of generations of bifurcations, since the relative size of red blood cells to vessel diameter increases, therefore the assumption made in 1D modelling that blood is a continuum and a Newtonian fluid fails.
Moreover, contrary to large arteries, fluid resistance dominates over wall compliance and fluid inertia in smaller vessels.
For these reasons, linear lumped parameter models are commonly employed to simulate the effect of peripheral vessels on pulse wave propagation in larger arteries~\cite{Alastruey2012a,westerhof1971artificial}.
Such models are obtained by averaging Eq.~\eqref{eq:1DAlastrueyAu} over the length of a vessel and considering some simplifying hypotheses, such as neglecting the convective term in the momentum equation.
A suitable lumped parameter model for our purposes is composed of one capacitor and two resistors \cite{Alastruey:2008a}.
In more detail, the first resistance  $Z$ corresponds to the characteristic aortic impedance.
This is connected in series with a parallel combination of the peripheral arterial resistance $R$ and compliance $C$.
We define $P_{out} $ as the pressure at which flow in the microcirculation is equal to zero.
	The resulting model is governed by
	\begin{equation}
		Q \Big(1 + \frac{Z}{R}\Big) + C\, Z \frac{\partial{Q}}{\partial{t}} = \frac{P - P_{out}}{R} + C \frac{\partial{P_e}}{\partial{t}}.
	\end{equation}
	\subsubsection{Numerical implementation}
	Numerous models have been proposed in recent years for the numerical discretisation of the resulting system of hyperbolic equations, relying on different discretisation approaches mainly based on the finite element method (FEM)~\cite{Formaggia_etal_2003,blanco-etal-14} or finite volume method~\cite{muller2014global}.
	A recent benchmark study on six commonly used numerical schemes for 1D blood flow modelling~\cite{boileau2015benchmark} showed a good agreement among these schemes and their capability to capture the main characteristics of flow and pressure waveforms in the large arteries.
In this study we used the solver~\href{http://haemod.uk/nektar}{Nektar1D} \cite{Alastruey2012a}, which is based on a FEM discretisation (discontinuous Galerkin) in space and finite difference (explicit second-order Adams--Bashforth) in time of Equations~\eqref{eq:1DAlastrueyAu} and~\eqref{eq:TubeLaw}.  
\section{Valve dynamics models}%
\label{sec:valves}
	The EM model of the heart is coupled to the arterial network via the cardiac valves.
	Flow through these valves is basically triggered by a positive difference pressure between the two compartments.
	However, several models can be considered for modelling the dynamics of these valves, ranging from simple diode-like models to more sophisticated approaches, that are able to contemplate pathological conditions such as a stenosis or regurgitation.
	In the literature there are three basic approaches, ranging from simple diode-like models to more complex valve models~\cite{quarteroni2017integrated}, simulating a simplified valve dynamics without explicitly modelling the valve leaflets.
	However, when it is necessary to capture more complex flow dynamics inside the ventricles that are induced by valve pathologies, it may be required to explicitly model the shape and dynamics of the leaflets.
	Several approaches have been proposed~\cite{votta2013toward}, but this is still a very open and demanding field of research.
	The simplest model to simulate the action of a valve corresponds to a diode associated with a resistance (characteristic of the valve), with instant opening and closing depending on the pressure gradient upstream and downstream the valve.
	This model has a physiological basis, since cardiac valves open and close very rapidly (in a few milliseconds).
	However, it does not take into account any dynamics on $Q$, and it is not possible to model pathological conditions.
	As a consequence, for this work we considered a more accurate model, namely a lumped parameter model of valve dynamics proposed in~\cite{Mynard_etal_2012}.
	Its peculiarity consists in a smooth opening and closing dynamics of the valves using a simple approach.
	In particular, this model is derived from the Bernoulli equation, after neglecting Poiseuille-type viscous losses (proportional to $Q$), since they are small:
	\begin{equation}
	\Delta P =  B\, Q\, |Q| + L \, \dot{Q}.
	\label{MynardValve}
	\end{equation}
	The Bernoulli resistance $B$ is associated with the pressure difference $\Delta P$ related to convective acceleration and dynamic pressure losses caused by the diverging flow field downstream to the vena contracta, and it reads
	\begin{equation}
	B = \frac{\rho}{2\, A_{\eff}^2},
	\label{eq:defB}
	\end{equation}
	where $\rho$ is the blood density (usually set equal to $1060\,\si{\kg\per\cubic\m}$) and $A_\eff$ is an effective cross-sectional area. The blood inertance $L$ is related to the pressure difference associated with blood acceleration and it is governed by the following equation:
	\begin{equation}
	L = \frac{\rho \, l_{\eff}}{A_{\eff}},
	\label{eq:defL}
	\end{equation}
	with $l_\eff$ an effective length, and it is here taken equal to the diameter of $A_\eff$.
	In order to consider valve dynamics, $A_\eff$ is a state variable and depends on an index of valve state $\xi \in [0,1]$ ($\xi = 0$: closed valve, $\xi = 1$: open valve) such that
	\begin{equation}
	A_\eff(t) = \big(  A_{\eff,\max}(t) - A_{\eff,\min}(t) \big)\xi (t) + A_{\eff,\min}(t).
	\end{equation}
	Note that, in order to avoid division by zero in Equations~\eqref{eq:defB} and~\eqref{eq:defL}, the valve state $\xi$ must be always strictly positive in practice.
	In order to consider pathological conditions like valve regurgitation, stenosis and varying annulus area $A_\mathrm{ann}$, the maximum $A_{\eff,\max}$ and minimum $A_{\eff,\min}$ effective areas can be expressed as
	\begin{equation}
	A_{\eff,\max}(t) = M_\mathrm{st}A_\mathrm{ann}(t), \quad A_{\eff,\min}(t) = M_\mathrm{rg}A_\mathrm{ann}(t).
	\end{equation}
	In this study, $A_\mathrm{ann}$ is considered constant. The coefficients $M_\mathrm{st}$ and $M_\mathrm{rg}$ range within $[0,1]$, with $M_\mathrm{rg}=0$ and $M_\mathrm{st}=1$ in case of a healthy valve, whereas $M_\mathrm{rg}>0$ indicates a regurgitant valve and $M_\mathrm{st}<1$ stands for a stenosed valve. Finally, $M_\mathrm{rg}=1$ corresponds to an absent valve (unobstructed conduit), whereas $M_\mathrm{st}=0$ stands for an atretic valve.\\
        The rate of valve opening and closing only depends on two factors, i.e., the instantaneous pressure difference across the valve $\Delta P$ and the current state $\xi$ of the valve. In particular, it is assumed that the valve opens when $\Delta P$ exceeds a threshold $\Delta P_\mathrm{open}$, and that the rate of opening is governed by
	\begin{equation}
	\frac{\der \xi}{\der t} = (1-\xi) K_\mathrm{vo}(\Delta P - \Delta P_\mathrm{open}),
	\end{equation}
	with $K_\mathrm{vo}$ a rate coefficient for valve opening (in \si{\per\pascal\per\second}). On the other hand,
	it is assumed that the valve closes when $\Delta P$ is lower than a threshold $\Delta P_\mathrm{close}$, and that the rate of closing is governed by
	\begin{equation}
	\frac{\der \xi}{\der t} = \xi\,  K_\mathrm{vc}(\Delta P - \Delta P_\mathrm{close}),
	\end{equation}
	where $K_{vc}$ is a rate coefficient for valve closing (in \si{\per\pascal\per\second}).
	Model parameters are depicted in~\Cref{tab:valve_input_params}.
		\begin{table}[htbp]
			\centering
			\footnotesize
			\begin{tabularx}{0.8\textwidth}{lrlL}
				\toprule
				Parameter & Value & Unit & Description \\
				\midrule
				\multicolumn{4}{l}{\emph{Aortic Valve }}\\
				$	K_\mathrm{vo}$           & 2   & \si{\per\pascal\per\second}    & rate coefficient for valve opening\\
				$\Delta P_\mathrm{open}$ & 0    & \si{\Pa} &  threshold pressure difference for valve opening\\
				$	K_\mathrm{vc}$ & 1.2 & \si{\per\pascal\per\second} & rate coefficient for valve closing \\
				$\Delta P_\mathrm{close}$ & 0    & \si{\Pa} &  threshold pressure difference for valve closing\\
				$M_\mathrm{st}$      & 1  & [-]    &  valve stenosis coefficient\\
				$M_\mathrm{rg}$      & 0  & [-]    &  valve regurgitation coefficient\\
				\midrule
				\multicolumn{4}{l}{\emph{Mitral Valve }}\\
				$	K_\mathrm{vo}$           & 0.2   & \si{\per\pascal\per\second}    & rate coefficient for valve opening\\
				$\Delta P_\mathrm{open}$ & 0    & \si{\Pa} &  threshold pressure difference for valve opening\\
				$	K_\mathrm{vc}$ & 0.2 & \si{\per\pascal\per\second} & rate coefficient for valve closing \\
				$\Delta P_\mathrm{close}$ & 0    & \si{\Pa} &  threshold pressure difference for valve closing\\
				$M_\mathrm{st}$      & 1  & [-]    &  valve stenosis coefficient\\
				$M_\mathrm{rg}$      & 0  & [-]    &  valve regurgitation coefficient\\
				\bottomrule
			\end{tabularx}
			\caption{Input parameters for the valve dynamics model.
				Adjusted to match patient-specific data.}%
			\label{tab:valve_input_params}
		\end{table}
\section{Coupling strategy}
\label{sec:coupling}
Various approaches can be envisaged to perform the coupling between a (3D) cardiac model and a circulatory model.
In practice, the problem consists in finding the new state of deformation $\vu^{n+1}$ as a function of the pressure $p^{n+1}$ in the cavity applied as a Neumann boundary condition at the cavitary surface.
This unknown pressure needs to be appropriately determined. %
Two configurations need to be captured by the model:
\begin{itemize}
 \item When all valves are closed, the cavity is in an isovolumetric state, i.e., the muscle enclosing the cavity may deform,
 but the volume has to remain constant.
 Consequently, during an isovolumetric phase the pressure $p^{n+1}$ in the cavity needs to adjust to the variation over time of active stresses, in order to maintain the cavitary volume constant.
 \item When one valve is open or regurging, there is a change in cavitary volume.
 In this configuration the pressure $p^{n+1}$ is regulated by the state of the circulatory system or of a connected cavity.
 Therefore, $p^{n+1}$ needs to be estimated by matching mechanical deformation
 and state of the system. In fact, pressure $p^{n+1}$ in the cardiovascular system
 depends on flow, which depends in turn by cardiac deformation.
 As a consequence, the heart and circulatory models are tightly bidirectionally coupled.
\end{itemize}
From a mathematical perspective, this coupling can be performed in two ways.
The simplest approach is to prescribe $p^{n+1}$ explicitly.
Then, the coupling is performed during the ejection phase by updating flow and flow rate
based on the current prediction on the change in the state of deformation
associated with the currently predicted pressure $p^{n+1}$.
In particular, the pressure boundary condition in each nonlinear solver step
is modified within each iteration ${\nu}$.
Note that the new prediction $p^{n+1}_{\nu+1}$ is then explicitly prescribed as a Neumann boundary condition.
Therefore, this explicit and partitioned approach may introduce some inaccuracies during ejection phases and can lead to instabilities during isovolumetric phases.
Such instabilities arise from the difficulty of estimating the change in pressure
necessary to keep the volume constant.
In fact, a knowledge on cavitary elastance is required to know the $p-V$ relation of the cavity at this given point in time, that is not available.
Therefore, iterative estimates are needed to gradually inflate or deflate
a cavity to its prescribed volume.
However, since the elastance properties of the cavities are highly nonlinear,
an overestimation may induce oscillations,
whereas an underestimation may induce very slow convergence and prohibitively large numbers of Newton iterations.
A more sophisticated approach consists in treating $p^{n+1}$ as an unknown.
Thus, an additional equation is required to close the system, leading to a saddle point formulation and a block system to be solved (monolithic approach).
\paragraph{Isovolumetric Phases} During the isovolumetric phases we can state that
\[
  \operatorname{V}^\mathrm{heart}_c(\vu)-V_0=0,
\]
i.e., the volume of each cavity ${V}^\mathrm{heart}_c(\vu)$ equals an initial volume $V_0$ and does not change during the isovolumetric phase.
Hence, in system~\eqref{eq:blockCVSystem} the matrix $\tensor{C}'(\Rvec{p}_k)=\tensor{0}$
and $V^{\mathrm{CS}}(\Rvec{p}_k)=V_{0}$.

\section{Finite Element Formulation}
\label{sec:FE}
\subsection{Variational Formulation}
\label{sec:fe_varf}
For the sake of simplicity of formulation, we neglect here the acceleration term in
Eq. (3) in this section and consider at the stationary counterpart
of the boundary value problem (\ref{eq:bvp}--\ref{eq:neumann}) with~\eqref{eq:VolumeEqu}.
Please refer to~\cite{augustin2021computationally} for further detail on the full nonlinear elastodynamics problem.
The boundary value problem (\ref{eq:bvp}--\ref{eq:neumann}) with~\eqref{eq:VolumeEqu} can be formally expressed in terms of the equations
\begin{align}
  \langle \mathcal{A}_0(\vu),\vv\rangle_{\Omega_0}
  - \langle \mathcal{F}_0(\vu,p_c),\vv\rangle_{\Omega_0} &=\vZero
  ,\label{eq:saddlePCVSystem1}\\
  \langle V_c^\mathrm{heart}(\vu),q\rangle_{\Omega_0}
  - \langle V_c^\mathrm{CS}(p_c), q\rangle_{\Omega_0}  &= 0,
  \label{eq:saddlePCVSystem2}
\end{align}
which is valid for all vector fields $\vv$ smooth enough and vanishing on
the Dirichlet boundary ${\Gamma}_{0,\mathrm{D}}$,
test functions $q$ that are $1$ for the cavity $c\in\{\mathrm{LV,RV,LA,RA}\}$ and $0$ otherwise, with
the duality pairing $\langle \cdot ,\cdot \rangle_{\Omega_0}$.
The first term on the left hand side of the variational
Eq.~\eqref{eq:saddlePCVSystem1} represents the
rate of internal mechanical work and is defined as
\begin{equation}\label{eq:lhs_stiffness}
    \langle \mathcal{A}_0(\vu) , \vv \rangle_{\Omega_0}
    :=\int_{\Omega_0} \tensor{S}(\vu) :
  \tensor{\Sigma}(\vu,\vv) \dd \vX,
\end{equation}
where $\tensor{S}$ denotes the second Piola--Kirchhoff stress tensor,  %
see~\eqref{eq:additiveSplit}, and $\tensor{\Sigma}(\vu,\vv)$ is the directional derivative of the Green--Lagrange strain tensor,  %
see~\cite{augustin2016anatomically, holzapfel2000nonlinear}.
The second term on the left hand side denotes the weak form of the contribution of pressure loads~\eqref{eq:saddlePCVSystem1} and is computed using~\eqref{eq:neumann}
\begin{equation} \label{eq:pressureLoads1}
  \langle \mathcal{F}_0(\vu,p_c),\vv \rangle_{\Omega_0} =
  -p_c \int\limits_{\Gamma_{0,\mathrm{N}}}J\,\tensor{F}^{-\top}(\vu)\,\normalout
  \cdot\vv \, \dsX.
\end{equation}
The first term of the coupling equation~\eqref{eq:saddlePCVSystem2} is
obtained from~\eqref{eq:volumeOmega} using Nanson's formula and
$\vx = \vX + \vu$ by
\begin{equation}\label{eq:VCAV}
  \langle V_c^\mathrm{heart}(\vu), q\rangle_{\Omega_0} =
  \frac{1}{3} \int\limits_{\Gamma_{0, \mathrm{N}}} (\vX + \vu)
  \cdot J \tensor F^{-\top} \normalout\, q \dsX.
\end{equation}
The second term of~\eqref{eq:saddlePCVSystem2} is computed using
the 1D blood flow model, see~\ref{sec:oneDmodelling},
for $c\in\{\mathrm{LV,RV,LA,RA}\}$.
\subsection{Consistent Linearisation}
A Newton--Raphson scheme is employed to solve the nonlinear variational
equations~\eqref{eq:saddlePCVSystem1}--\eqref{eq:saddlePCVSystem2},
with a finite element approach %
see~\cite{deuflhard2011newton}.
Let us consider a nonlinear and continuously differentiable operator \(F\colon X\to Y\). Then,
a solution to \(F(x)=0\) can be approximated by
\begin{align*}
  x^{k+1} &= x^{k} + \Delta x, \\
  \left.\frac{\partial F}{\partial x}\right|_{x = x^k} \Delta x &= -F(x^k),
\end{align*}
which is evaluated until convergence.
In our approach, we assume \(X = \left[H^1(\Omega_0,\Gamma_{0,\mathrm{D}})\right]^3 \times \mathbb R\), \(Y = \mathbbm R^2\),
\(\Delta x = {(\Delta \vu, \Delta p_c)}^\top\), \(x^k = {(\vu^k, p_c^k)}^\top\),
and \(F = {(R_{\vu},R_\mathrm{p})}^\top\).
Thus, we obtain the following linearised saddle-point problem for each
\((\vec u^k, p_c^k) \in \left[H^1(\Omega_0,\Gamma_{0,\mathrm{D}})\right]^3
\times \mathbb R\),
find \((\Delta \vec u, \Delta p_c) \in \left[H^1_0(\Omega_0)\right]^3
\times \mathbb R\) such that
\begin{align}
  \label{eq:saddle_point_nl:1}
  \langle \Delta\vu,A_0'(\vu^k)\, \vv\rangle_{\Omega_0}
  + \langle \Delta\vu,\mathcal{F}_0'(\vu^k, p_c^k)\,\vv\rangle_{\Omega_0}&\nonumber
  \\+ \langle \Delta p_c,\mathcal{F}_0'(\vu^k, p_c^k)\,\vv\rangle_{\Omega_0}
  &= -\langle R_{\vu}(\vu^k, p_c^k), \vv\rangle_{\Omega_0},\\
  \label{eq:saddle_point_nl:2}
  \langle \Delta\vu,V_c^\mathrm{heart}(\vu^k)\, q\rangle_{\Omega_0}
  - \langle \Delta p_c, V_c^\mathrm{CS}(p_c^k)\, q \rangle_{\Omega_0} &=
  -\langle R_\mathrm{p}(\vu^k, p_c^k), q\rangle_{\Omega_0},
\end{align}
with the updates
\begin{align}
  \vu^{k+1}&=\vu^k+\Delta\vu,\\
  p_c^{k+1}&=p_c^k+\Delta p_c.
\end{align}
The G\^ateaux derivative of~\eqref{eq:saddlePCVSystem1} with respect to the
displacement change update $\Delta \vu$ yields the first term in~\eqref{eq:saddle_point_nl:1}
\begin{align}
  \langle \Delta\vu,A_0'(\vu^k)\, \vv\rangle_{\Omega_0}
  :&= \left.D_{\Delta\vu} \langle \mathcal{A}_0(\vu),
        \vv\rangle_{\Omega_0}\right|_{\vu =\vu^k} \nonumber\\
  &= \int\limits_{\Omega_0} \tensor{S}_k: \tensor \Sigma(\Delta \vu, \vv)\dX
    + \int\limits_{\Omega_0}\tensor \Sigma(\vu^k, \Delta \vu)
    : \mathbbm C_k : \tensor \Sigma(\vu^k, \vv)\dX,
    \label{eq:linearizedSaddlePCV1}
\end{align}
and the second term in~\eqref{eq:saddle_point_nl:1}, that is given by
\begin{align}
  \langle \Delta\vu,\mathcal{F}_0'(\vu^k, p_c^k)\,\vv\rangle_{\Omega_0}
  :&= \left.D_{\Delta\vu} \langle \mathcal{F}_0(\vu,p_c),
    \vv\rangle_{\Omega_0}\right|_{\vec u =\vec u^k, p_c = p_c^k}\nonumber\\
  &= p^k_c \int_{\Gamma_{0,\mathrm{N}}}J_k\tensor{F}_k^{-\top}\Grad^\top\!\Delta\vu\,
  \tensor{F}_k^{-\top}\normalout\cdot\vv\, \dsX \nonumber \\
  &-p^k_c \int_{\Gamma_{0,\mathrm{N}}}J_k(\tensor{F}_k^{-\top}:\Grad\Delta\vu)\,
  \tensor{F}_k^{-\top}\normalout\cdot\vv\, \dsX,
    \label{eq:linearizedSaddlePCV2}
\end{align}
where we have defined, for the sake of simplicity:
\begin{align*}
  \tensor F_k &:= \tensor F(\vu^k),\
  J_k := \det(\tensor F^k),\
  \tensor{S}_k:= \left.\tensor S\right|_{\vu =\vu^k},\
  \mathbbm{C}_k:= \left.\mathbbm{C}\right|_{\vu =\vu^k}.
\end{align*}
The third term in~\eqref{eq:saddle_point_nl:1} is defined using the G\^ateaux derivative of~\eqref{eq:saddlePCVSystem1} with respect to the
pressure change update $\Delta p_c$:
\begin{align}
  \langle \Delta p_c,\mathcal{F}_0'(\vu^k, p_c^k)\,\vv\rangle_{\Omega_0}
:&=\left.D_{\Delta p_c} \langle \mathcal{F}_0(\vu,p_c),
        \vv\rangle_{\Omega_0}\right|_{\vec u =\vec u^k, p_c = p_c^k} \nonumber\\
        &=-\Delta p_c  \int_{\Gamma_{0,\mathrm{N}}}J_k\,\tensor{F}_k^{-\top}\,\normalout\cdot\vv \, \dd s_{\vX}.
    \label{eq:linearizedSaddlePCV3}
\end{align}
The right hand side of~\eqref{eq:saddle_point_nl:1} depends on the residual $R_{\vu}$, and it is defined as
\begin{equation}\label{eq:upper_rhs_variational}
\langle R_{\vu}(\vu^k, p_c^k), \vv\rangle_{\Omega_0}
  :=\langle A_0(\vu^k),\vv\rangle_{\Omega_0}
  - \langle \vec{\mathcal{F}}_0(\vu^k, p_c^k),\vv \rangle_{\Omega_0}.
\end{equation}
From~\eqref{eq:VCAV},
using the known relations~\cite{holzapfel2000nonlinear},
\begin{align*}
  \frac{\partial J}{\partial \tensor F} : \Grad \Delta \vu
    &= J \tensor F^{-\top} : \Grad \Delta \vu \\
  \frac{\partial \tensor F^{-\top}}{\partial \tensor F} : \Grad \Delta \vu
  &= -\tensor F^{-\top}{(\Grad \Delta \vu)}^\top \tensor F^{-\top}
\end{align*}
we can define the first term of~\eqref{eq:saddle_point_nl:2} as the
G\^ateaux derivative with respect to the update $\Delta\vu$
\begin{align}
  \langle \Delta\vu,V_c^\mathrm{heart}(\vu^k)\, q\rangle_{\Omega_0}
  :&=\left.D_{\Delta \vu}
  \langle V_c^\mathrm{heart}(\vu),q\rangle_{\Omega_0}
    \right|_{\vu =\vu^k} \nonumber \\
    &= D_{\Delta \vu} \frac{1}{3}\int_{\Gamma_{0,\mathrm{N}}}
    \left(\vX + \vu^k\right) \cdot J_k \tensor F_k^{-\top}
    \normalout\,q \dsX \nonumber \\
    &= \frac{1}{3} \int_{\Gamma_{0,\mathrm{N}}} J_k (\tensor F_k^{-\top}
  : \Grad \Delta \vu) \vx \cdot \tensor F_k^{-\top}
    \normalout\,q \dsX \nonumber \\
    &\quad - \frac{1}{3} \int_{\Gamma_{0,\mathrm{N}}} J_k \vx \cdot
    \tensor F_k^{-\top}{(\Grad \Delta \vu)}^\top
    \tensor F_k^{-\top} \normalout\,q \dsX \nonumber \\
    &\quad + \frac{1}{3} \int_{\Gamma_{0,\mathrm{N}}} J_k \Delta \vu \cdot
    \tensor F_k^{-\top}\normalout\,q \dsX,  \label{eq:var_volume_cav}
\end{align}
where $q$ is a test function with value $1$ for the surface $\Gamma_{0,c}$ of cavity $c$,
and $0$ otherwise.

The second term of~\eqref{eq:saddle_point_nl:2} is defined as the
numerical derivative
\begin{align}
  \langle \Delta p_c, V_c^\mathrm{CS}(p_c^k)\, q \rangle_{\Omega_0}
  :&= \left.D_{\Delta p_c} \langle V_c^\mathrm{CS}(p_c)\, q \rangle_{\Omega_0}
    \right|_{p_c = p_c^k} \nonumber \\
    &= \frac{1}{\epsilon}\left(V_c^\mathrm{CS}(p_c^k + \epsilon) -
    V_c^\mathrm{CS}(p_c^k)\right)q, \label{eq:var_compliance}
\end{align}
with $\epsilon=p_c^k\sqrt{\epsilon_\mathrm{m}}$ chosen according
to~\cite[Chapter 5.7]{Press2007}, and
$\epsilon_\mathrm{m}=2^{-52}\approx2.2*10^{-16}$ is the machine accuracy.

Finally, the right hand side of~\eqref{eq:saddle_point_nl:2}, is defined based on the residual $R_\mathrm{p}$ and reads
\begin{equation} \label{eq:lower_rhs_variational}
\langle R_\mathrm{p}(\vu^k, p_c^k), q\rangle_{\Omega_0}
  := \langle V_c^\mathrm{heart}(\vu),q\rangle_{\Omega_0}
  - \langle V_c^\mathrm{CS}(p_c), q\rangle_{\Omega_0}.
\end{equation}

\paragraph{Assembling of the block matrices}%
\label{par:blockMatrices}
The finite element method (FEM) relies on the definition of an admissible decomposition
of the computational domain $\Omega \subset \mathbbm{R}^3$
into $M$ tetrahedral elements $\tau_j$ and the introduction of a conformal finite
element space
\begin{align*}
X_h \subset H^1(\Omega_0), \ N = \text{dim}\, X_h
\end{align*}
of piecewise polynomial continuous basis functions $\varphi_i$.
After linearising the variational problem~\eqref{eq:saddle_point_nl:1}--\eqref{eq:saddle_point_nl:2}
and considering a Galerkin finite element discretisation, we obtain a
block system to be solved.
The problem then reads:
find $\delta\Rvec{u}\in\mathbbm{R}^{3N}$ and
$\delta \Rvec{p}_c\in\mathbbm{R}^{N_\mathrm{cav}}$ such that
\begin{align}\label{eq:blockCVSystem2}
  \begin{pmatrix}
    (\tensor{A}'-\tensor{M}')(\Rvec{u}^k, \Rvec{p}_c^k)  & \tensor{B}'_\mathrm{p}(\Rvec{u}^k) \\
        \tensor{B}'_{\vu}(\Rvec{u}^k) & \tensor{C}'(\Rvec{p}_c^k) \\
 \end{pmatrix}
\begin{pmatrix}
  \delta\Rvec{u}\\ \delta \Rvec{p}_c
\end{pmatrix}
&=-
\begin{pmatrix}
  \Rvec{A}(\Rvec{u}^k)-\Rvec{B}_\mathrm{p}(\Rvec{u}^k,\Rvec{p}_c^k) \\
  \Rvec{V}_c^\mathrm{heart}(\Rvec{u}^k) -\Rvec{V}_c^\mathrm{CS}(\Rvec{p}_c^k)
\end{pmatrix},\\
    \Rvec{u}^{k+1}   &= \Rvec{u}^{k} + \delta\Rvec{u},\\
    \Rvec{p}_c^{k+1} &= \Rvec{p}_c^k + \delta\Rvec{p}_c
\end{align}
with the solution vectors $\Rvec{u}^k\in\mathbbm{R}^{3N}$ and
$\Rvec{p}_c^k\in\mathbbm{R}^{N_\mathrm{cav}}$ at the $k$-th Newton
step.
The tangent stiffness matrix $\tensor{A}'\in\mathbbm{R}^{3N\times 3N}$
is computed from~\eqref{eq:linearizedSaddlePCV1} as follows:
\begin{equation}\label{eq:lhs_assemblingv}
  \tensor{A}'(\Rvec{u}^k)[j,i] :=
    \langle \boldsymbol{\varphi}_i, \mathcal{A}_0'(\vu^k)\,
            \boldsymbol{\varphi}_j \rangle_{\Omega_0}.
  \end{equation}
The mass matrix $\tensor{M}'\in\mathbbm{R}^{3N\times 3N}$
is calculated from~\eqref{eq:linearizedSaddlePCV2} and reads
\begin{equation}\label{eq:lhs_mass}
  \tensor{M}'(\Rvec{u}^k, \Rvec{p}_c^k)[j,i] :=
    \langle \boldsymbol{\varphi}_i, \mathcal{F}_0'(\vu^k, p_c^k)\,
            \boldsymbol{\varphi}_j \rangle_{\Omega_0}.
\end{equation}
The off-diagonal matrices
$\tensor{B}'_{\vu}\in\mathbbm{R}^{3N\times N_\mathrm{cav}}$ and
$\tensor{B}'_{\mathrm{p}}\in\mathbbm{R}^{N_\mathrm{cav}\times 3N}$
in~\eqref{eq:blockCVSystem2} are assembled using~\eqref{eq:var_volume_cav}
\begin{equation}\label{eq:Bu_assembling}
  \tensor{B}'_{\vu}(\Rvec{u}^k,\Rvec{p}_c^k)[i,j]
  = \langle \underline{\varphi}_j,
  V_c^\mathrm{heart}(\vu^k)\hat{\varphi_i}\rangle_{\Omega_{0}},\quad
            i=1,\dots,N_\mathrm{cav}
\end{equation}
and~\eqref{eq:linearizedSaddlePCV3}
\begin{equation}\label{eq:Bp_assembling}
  \tensor{B}'_{\mathrm{p}}(\Rvec{u}^k,\Rvec{p}_c^k)[i,j]
  = \langle \hat{\varphi_j},\mathcal{F}_0'(\vu^k, p_c^k)
  \boldsymbol{\varphi}_i\rangle_{\Omega_0}, \quad
  j=1,\dots,N_\mathrm{cav},
\end{equation}
where the constant shape functions have value $\hat{\varphi_j}=1$ if $\tau_j\in\Gamma_{0,c}$
and $\hat{\varphi_j}=0$ if $\tau_j\notin\Gamma_{0,c}$
for $c\in\{\mathrm{LV,RV,LA,RA}\}$.

Considering the approach described in~\cite[Sect.~4.2]{rumpel2003volume}, this
assembling procedure can be simplified for closed cavities such that
\[
  \tensor{B}'_{\mathrm{p}}(\Rvec{u}^k,\Rvec{p}_c^k) =
  {\left[\tensor{B}'_{\vu}(\Rvec{u}^k,\Rvec{p}_c^k)\right]}^\top.
\]
Then, the circulatory compliance matrix
$\tensor{C}'(\Rvec{p}_c^k)\in\mathbbm{R}^{N_\mathrm{cav}\times N_\mathrm{cav}}$
is computed from~\eqref{eq:var_compliance} as
\begin{equation}\label{eq:Cp_assembling}
  \tensor{C}'(\Rvec{p}_c^k)[i,j]=
  \langle \hat{\varphi_j}, V_c^\mathrm{CS}(p_c^k)\,
    \hat{\varphi_i}\rangle_{\Omega_0}, \quad
    i,j=1,\dots,N_\mathrm{cav},
\end{equation}
with the constant shape function $\hat{\varphi_j}=1$ if $\tau_j\in\Gamma_{0,c}$
and $\hat{\varphi_j}=0$ if $\tau_j\notin\Gamma_{0,c}$
for $c\in\{\mathrm{\LV,\RV,\LA,\RA}\}$, leading to a diagonal matrix.
The terms on the right hand side $\Rvec{A}\in\mathbbm{R}^{3N}$,
$\Rvec{B}_\mathrm{p}\in\mathbbm{R}^{3N}$ are constructed
using~\eqref{eq:upper_rhs_variational} and read
\begin{equation}\label{eq:residual_A}
\Rvec{A}(\Rvec{u}^k)[i] := \langle
\mathcal{A}_0(\vu^k),\boldsymbol{\varphi}_i \rangle_{\Omega_0}
\end{equation}
and
\begin{equation}\label{eq:residual_B}
\Rvec{B}_\mathrm{p}(\Rvec{u}^k, \Rvec{p}_c^k)[i] := \langle
\mathcal{F}_0(\vu^k, p_c^k),\boldsymbol{\varphi}_i \rangle_{\Omega_0},
\end{equation}
see also~\cite{augustin2016anatomically, holzapfel2000nonlinear}.
Finally, the terms $\Rvec{V}_c^\mathrm{heart}\in\mathbbm{R}^{N_\mathrm{cav}}$ and
$\Rvec{V}_c^\mathrm{CS}\in\mathbbm{R}^{N_\mathrm{cav}}$ in lower right hand side in~\eqref{eq:blockCVSystem2} are assembled from~\eqref{eq:lower_rhs_variational} and are given by:
\begin{equation}\label{eq:VcPDE_assembling}
  \Rvec{V}_c^\mathrm{heart}(\Rvec{u}^k)[i]
  =\langle V_c^\mathrm{heart}(\vu),\hat{\varphi_i}\rangle_{\Omega_0},\quad
    i=1,\dots,N_\mathrm{cav},
\end{equation}
and
\begin{equation}\label{eq:VcCS_assembling}
  \Rvec{V}_c^\mathrm{CS}(\Rvec{p}_c^k)[i]
  =\langle V_c^\mathrm{CS}(p_c), \hat{\varphi_i}\rangle_{\Omega_0},\quad
    i=1,\dots,N_\mathrm{cav}.
\end{equation}
\section{Direct Schur Complement Solver for a Small Number of Constraints}
\label{sec:SchurComplement}
Let us consider the block system $\tensor{K}\in\mathbbm{R}^{n\times n}$, $\tensor{C}\in\mathbbm{R}^{m\times m}$
\vspace{-0.05cm}
\begin{align*}
  \begin{pmatrix}
    \tensor{K}&\tensor{A}\\
    \tensor{B}&\tensor{C}
  \end{pmatrix}
  \begin{pmatrix}
    \Rvec{u}\\ \Rvec{p}
  \end{pmatrix}=-
  \begin{pmatrix}
    \Rvec{f}\\ \Rvec{g}
  \end{pmatrix}
\end{align*}
with
\vspace{-0.05cm}
\begin{align*}
    \tensor{A}=\begin{pmatrix} \Rvec{a}_1 & \aug & \cdots & \aug &\Rvec{a}_m\end{pmatrix}\in\mathbbm{R}^{n\times m}, \quad
  \tensor{B}=\begin{pmatrix} \Rvec{b}_1 & \aug & \cdots & \aug & \Rvec{b}_m\end{pmatrix}^\top\in\mathbbm{R}^{m\times n}.
\end{align*}
We can write the Schur complement system as
\vspace{-0.05cm}
\begin{align*}
  (\tensor{B}\tensor{K}^{-1}\tensor{A}-\tensor{C})\Rvec{p}&=\Rvec{g}-\tensor{B}\tensor{K}^{-1}\Rvec{f}\\
\Rvec{u}&=\tensor{K}^{-1}\Rvec{f}-\tensor{K}^{-1}\tensor{A}\Rvec{p}.
\end{align*}
With
\vspace{-0.05cm}
\begin{equation} \label{}
  \Rvec{r}=\tensor{K}^{-1}\Rvec{f},\quad
  \tensor{S}=\tensor{K}^{-1}\Rvec{A}=
  \begin{pmatrix}\Rvec{s}_1 & \aug & \cdots & \aug & \Rvec{s}_m\end{pmatrix}\in\mathbbm{R}^{n\times m},\
  \Rvec{s}_i=\tensor{K}^{-1}\Rvec{a}_i,\ i=1,\dots,m
\end{equation}
we obtain
\vspace{-0.05cm}
\begin{align}
  (\tensor{B}\tensor{S}-\tensor{C})\Rvec{p}&=\Rvec{g}-\tensor{B}\Rvec{r}\nonumber\\
\Rvec{u}&=\Rvec{r}-\tensor{S}\Rvec{p}.\label{eq:SchurSolve}
\end{align}
The realisation of~\eqref{eq:SchurSolve} involves $m+1$ solves and the inversion of an $m\times m$ matrix. As $m$ is generally small, this can be done symbolically.
\vspace{-0.05cm}
\begin{align*}
    {[\tensor{B}\tensor{S}]}_{ij}=\Rvec{b}_i\cdot\Rvec{s}_j,\ \mbox{for}\ i,j=1,\dots,m.
\end{align*}
\end{appendix}
\end{document}